\documentclass[11pt, a4paper]{article}
\usepackage[cp1251]{inputenc}
\usepackage{amsmath} \usepackage{euscript}
\usepackage{amssymb}
\begin{document}

\newcounter{bnomer} \newcounter{snomer}
\newcounter{bsnomer}
\setcounter{bnomer}{-1}
\renewcommand{\thesnomer}{\thebnomer.\arabic{snomer}}
\renewcommand{\thebsnomer}{\thebnomer.\arabic{bsnomer}}
\renewcommand{\refname}{\begin{center}\large{\textbf{References}}\end{center}}

\newcommand{\sect}[1]{%
\setcounter{snomer}{0}\setcounter{bsnomer}{0}
\refstepcounter{bnomer}
\par\bigskip\begin{center}\large{\textbf{\arabic{bnomer}. {#1}}}\end{center}}
\newcommand{\sst}{%
\refstepcounter{bsnomer}
\par\bigskip\textbf{\arabic{bnomer}.\arabic{bsnomer}. }}
\newcommand{\defi}[1]{%
\refstepcounter{snomer}
\par\textbf{Definition \arabic{bnomer}.\arabic{snomer}. }{#1}\par}
\newcommand{\theo}[2]{%
\refstepcounter{snomer}
\par\textbf{Theorem \arabic{bnomer}.\arabic{snomer}. }{#2} {\emph{#1}}\hspace{\fill}$\square$\par}
\newcommand{\mtheo}[1]{%
\refstepcounter{snomer}
\par\textbf{Theorem \arabic{bnomer}.\arabic{snomer}. }{\emph{#1}}\par}
\newcommand{\theobp}[2]{%
\refstepcounter{snomer}
\par\textbf{Theorem \arabic{bnomer}.\arabic{snomer}. }{#2} {\emph{#1}}\par}
\newcommand{\theop}[2]{%
\refstepcounter{snomer}
\par\textbf{Theorem \arabic{bnomer}.\arabic{snomer}. }{\emph{#1}}
\par\textbf{Proof}. {#2}\hspace{\fill}$\square$\par}
\newcommand{\theosp}[2]{%
\refstepcounter{snomer}
\par\textbf{Theorem \arabic{bnomer}.\arabic{snomer}. }{\emph{#1}}
\par\textbf{Sketch of the proof}. {#2}\hspace{\fill}$\square$\par}
\newcommand{\exam}[1]{%
\refstepcounter{snomer}
\par\textbf{Example \arabic{bnomer}.\arabic{snomer}. }{#1}\par}
\newcommand{\deno}[1]{%
\refstepcounter{snomer}
\par\textbf{Definition \arabic{bnomer}.\arabic{snomer}. }{#1}\par}
\newcommand{\post}[1]{%
\refstepcounter{snomer}
\par\textbf{Proposition \arabic{bnomer}.\arabic{snomer}. }{#1}\hspace{\fill}$\square$\par}
\newcommand{\postp}[2]{%
\refstepcounter{snomer}
\par\textbf{Proposition \arabic{bnomer}.\arabic{snomer}. }{\emph{#1}}
\par\textbf{Proof}. {#2} $\square$\par}
\newcommand{\lemm}[1]{%
\refstepcounter{snomer}
\par\textbf{Lemma \arabic{bnomer}.\arabic{snomer}. }{\emph{#1}}\hspace{\fill}$\square$\par}
\newcommand{\lemmp}[2]{%
\refstepcounter{snomer}
\par\textbf{Lemma \arabic{bnomer}.\arabic{snomer}. }{\emph{#1}}
\par\textbf{Proof}. {#2}\hspace{\fill}$\square$\par}
\newcommand{\coro}[1]{%
\refstepcounter{snomer}
\par\textbf{Corollary \arabic{bnomer}.\arabic{snomer}. }{\emph{#1}}\hspace{\fill}$\square$\par}
\newcommand{\mcoro}[1]{%
\refstepcounter{snomer}
\par\textbf{Corollary \arabic{bnomer}.\arabic{snomer}. }{\emph{#1}}\par}
\newcommand{\corop}[2]{%
\refstepcounter{snomer}
\par\textbf{Corollary \arabic{bnomer}.\arabic{snomer}. }{\emph{#1}}
\par\textbf{Proof}. {#2}\hspace{\fill}$\square$\par}
\newcommand{\nota}[1]{%
\refstepcounter{snomer}
\par\textbf{Remark \arabic{bnomer}.\arabic{snomer}. }{#1}\par}
\newcommand{\propp}[2]{%
\refstepcounter{snomer}
\par\textbf{Proposition \arabic{bnomer}.\arabic{snomer}. }{\emph{#1}}
\par\textbf{Proof}. {#2}\hspace{\fill}$\square$\par}

\newcommand{\Ind}[3]{%
\mathrm{Ind}_{#1}^{#2}{#3}}
\newcommand{\Res}[3]{%
\mathrm{Res}_{#1}^{#2}{#3}}
\newcommand{\epsi}{\varepsilon}
\newcommand{\Supp}[1]{%
\mathrm{Supp}(#1)}

\newcommand{\reg}{\mathrm{reg}}
\newcommand{\sreg}{\mathrm{sreg}}
\newcommand{\codim}{\mathrm{codim}\,}
\newcommand{\chara}{\mathrm{char}\,}
\newcommand{\rk}{\mathrm{rk}\,}
\newcommand{\id}{\mathrm{id}}
\newcommand{\col}{\mathrm{col}}
\newcommand{\row}{\mathrm{row}}
\newcommand{\pho}{\hphantom{\quad}\vphantom{\mid}}
\newcommand{\wt}{\widetilde}
\newcommand{\wh}{\widehat}
\newcommand{\ad}[1]{\mathrm{ad}_{#1}}

\newcommand{\vfi}{\varphi}
\newcommand{\teta}{\vartheta}
\newcommand{\lee}{\leqslant}
\newcommand{\gee}{\geqslant}
\newcommand{\Fp}{\mathbb{F}}
\newcommand{\Rp}{\mathbb{R}}
\newcommand{\Zp}{\mathbb{Z}}
\newcommand{\Cp}{\mathbb{C}}
\newcommand{\dl}{\delta}
\newcommand{\ut}{\mathfrak{u}}
\newcommand{\at}{\mathfrak{a}}
\newcommand{\rt}{\mathfrak{r}}
\newcommand{\rad}{\mathfrak{rad}}
\newcommand{\bt}{\mathfrak{b}}
\newcommand{\gt}{\mathfrak{g}}
\newcommand{\vt}{\mathfrak{v}}
\newcommand{\pt}{\mathfrak{p}}
\newcommand{\Po}{\EuScript{P}}
\newcommand{\Uo}{\EuScript{U}}
\newcommand{\Fo}{\EuScript{F}}
\newcommand{\Mo}{\mathcal{M}}
\newcommand{\Ro}{\mathcal{R}}
\newcommand{\Co}{\mathcal{C}}
\newcommand{\Lo}{\mathcal{L}}
\newcommand{\Ou}{\mathcal{O}}
\newcommand{\Au}{\mathcal{A}}
\newcommand{\Bu}{\mathcal{B}}
\newcommand{\Sy}{\mathcal{Z}}
\newcommand{\Sb}{\mathcal{F}}
\newcommand{\Gr}{\mathcal{G}}

\hyphenation{Pro-po-si-tion}

\author{Mikhail V. Ignatyev\thanks{Samara state university,
Department of algebra and geometry, 443011, ak. Pavlova, 1, Samara,
Russia, \texttt{mihail.ignatev@gmail.com}}}
\date{}
\title{Orthogonal subsets of root systems and the orbit method} \maketitle

\sect{Introduction and statement of the main result}


\sect{Introduction}

\sst The main tool in studying irreducible complex representations
of finite unipotent groups is the \emph{orbit method}. It was
created by A.A. Kirillov for nilpotent Lie groups over $\Rp$
\cite{Kirillov1}, \cite{Kirillov2}, and~then adapted by D.~Kazhdan
for finite groups~\cite{Kazhdan} (see also \cite{Kirillov3} and the
paper \cite{BoyDr}, where the theory of $\ell$-adic sheaves for
unipotent groups is explained). Here we consider the groups $U(q)$
and $U$, the maximal unipotent subgroups of Chevalley groups over a
finite field~$\Fp_q$ and its algebraic closure respectively.

The orbit method establishes a bijection between the set of
equivalence classes of~irreducible representations of $U(q)$ and the
set of~orbits of the \emph{coadjoint representation} of $U(q)$.
Further, a lot of questions about {re\-pre\-sen\-ta\-tions} can be
interpreted in terms of~orbits. Note that the~problem of complete
description of orbits remains unsolved and seems to be \emph{very}
difficult. On~the~other hand, a~lot of information about some
special types of orbits, {re\-pre\-sen\-ta\-tions} and characters is
known.

In particular, a description of \emph{regular} orbits (i.e., orbits
of maximal {di\-men\-sion}) of the group $\mathrm{UT}_n$ of all
unipotent triangular matrices of size $n\times n$ is known
\cite{Kirillov2}. \emph{Subregular} orbits (i.e., orbits of~second
maximal {di\-men\-sion}) and corresponding characters\footnote{The
description of the irreducible character corresponding to a given
orbit is itself a~non-trivial problem, see, e.g., the papers of
C.A.M. Andr\`e and A. Neto \cite{Andre}, \cite{AndreNeto} for the
description of the so-called \emph{super}-\emph{characters}. In this
paper, we concentrate on orbits, not on characters.} were described
in \cite{IgnatevPanov} and \cite{Ignatev1}. As a generalization,
A.N. Panov considered orbits of the group $\mathrm{UT}_n$ associated
with involutions in the symmetric group. In \cite{Panov}, he
obtained a formula for the dimension of such an orbit.

It's well-known that the group $\mathrm{UT}_n$ corresponds to the
root system of type $A_{n-1}$. In order to generalize the~results of
A. N. Panov, we introduced the concept of orbits associated with
orthogonal subsets of root systems. In the paper \cite{Ignatev3}, we
studied these orbits for~the~case of classical root systems.
For~orthogonal subsets of special kind of the~root systems of types
$B_n$~and~$D_n$, we also obtained a~formula involving the
corresponding irreducible characters, see \cite[Theorem
3.8]{Ignatev2}).

The main goal of this paper is to generalize the results of
\cite{Ignatev3} to the general case of an \emph{arbitrary} root
system, not only the classical one. The structure of the paper ia as
follows. In the remainder of this Section, we give necessary
definitions and formulate the main result (see Theorem \ref{mtheo}).
In Section~1, we prove some preliminary technical Lemmas and
consider some important examples. In Section~2, we~prove the Main
Theorem for simply laced root systems
(see~Propositions~\ref{prop_dim_s_laced}~and~\ref{prop_xi_s_laced}).
In Section~3, we prove the Main Theorem for multiply laced root
systems.

The author is sincerely grateful to his scientific advisor professor
A.N.~Panov for constant attention to this work.

\sst In this Subsection, we shall briefly recall some basic facts
concerning Chevalley groups over finite fields. We also give some
definitions, which are needed to formulate the main result.

Let $\Phi$ be a reduced root system, $\Delta\subset\Phi$ a subset of
fundamental roots, $\Phi^+$~and~$\Phi^-$ the corresponding subsets
of positive and negative roots {respec\-tively}
(see~\cite{Bourbaki}). As usual, we denote by~${W=W(\Phi)}$ the~Weyl
group of the root system $\Phi$. Let $r_{\alpha}\in W$~be the
reflection on the hyperplane orthogonal to~a~given
root~$\alpha\in\Phi$.

Let $p$ be a prime, $\Fp_q$ the field with $q=p^r$ elements for some
${r\gee1}$, $k=\overline{\Fp}_q$ its algebraic closure. Let
$G(q)=G_{\mathrm{sc}}(\Phi, \Fp_q)$ (resp. $G=G_{\mathrm{sc}}(\Phi,
k)$) be~the~simply connected \emph{Chevalley group} over the field
$\Fp_q$ (resp. over $k$) with the root {sys\-tem}~$\Phi$
(see~the~classical book~\cite{Steinberg1} for precise definitions;
see~also~\cite{Srinivasan}). Recall that there exists a so-called
\emph{Chevalley basis} of the Lie algebra~$\gt$ of the group $G$. In
particular, this basis contains the root vectors $\{e_{\alpha}$,
$\alpha\in\Phi^+\}$ satisfying
${[e_{\alpha},e_{\beta}]=N_{\alpha\beta}e_{\alpha+\beta}}$, where
$N_{\alpha\beta}$ are the so-called \emph{Chevalley structure
constants} (here we set $N_{\alpha\beta}=0$
if~$\alpha+\beta\notin\Phi$).

The subspace $\ut=\sum_{\alpha\in\Phi^+}ke_{\alpha}$ is a nilpotent
Lie subalgebra of $\gt$. We assume from now on that $p$ is not less
than the Coxeter number of the root system $\Phi$. This implies
$[x_1,[x_2,[\ldots,[x_{p-1}, x_p]\ldots]]=0$ for all $x_i\in\ut$, so
the orbit method applies \cite[Theorem 2.2 and \S3.3]{BoyDr}.

Since $p$ is sufficiently large, the~exponential map
$\exp\colon\ut\to G$ is well-defined. Its image $U$ is a maximal
unipotent subgroup of $G$ and the~map $\exp\colon\ut\to U$ is a
bijection. Further, $U$ is generated as a subgroup of $G$ by all
root subgroups corresponding to positive roots from $\Phi$, and
$\ut$ is the Lie algebra of the group $U$.

Thus, the group $U$ acts on $\ut$ via the adjoint representation.
The dual {re\-pre\-sen\-ta\-tion} of $U$ in~the space $\ut^*$ of all
$k$-linear functions on $\ut$ is called \emph{coadjoint}. One can
see that the~coadjoint action has the form
\begin{equation*}
\exp(y).f(x)=f(\exp\ad{-y}x),\quad x,y\in\ut,f\in\ut^*.
\end{equation*}
Here $\ad{y}x=[y, x]$; since $\ad{y}$ is a nilpotent linear operator
on $\ut$, the map
$\exp\ad{y}=\sum_{i=0}^{\infty}\ad{y}^i/i!\colon\ut\to U$
is~well-defined.

One can define the algebra $\ut(q)\subset\gt(q)$, the group $U(q)$
and its coadjoint representation in~the~space
${\ut^*(q)=(\ut(q))^*}$ by the similar way. Let us fix an embedding
$\Fp_q\subset k$. Then $\ut(q)$ can be canonically embedded  in
$\ut$. In the paper we concentrate\footnote{The set of irreducible
representations of the group $U(q)$ is in bijection with the set of
coadjoint orbits of $U(q)$, but a lot of questions about
representations can be interpreted in terms of orbits of the group
$U$ \cite{Kazhdan}.} on orbits of elements from $\ut(q)$ under the
coadjoint action of the group $U$, not of the group $U(q)$.

Now we shall give the main definition. Let $D$ be a subset of
$\Phi^+$ consisting of pairwise orthogonal roots, then $D$ is called
\emph{orthogonal}. Let $\xi=(\xi_{\beta})_{\beta\in D}$ be a set of
non-zero scalars from $k$. Denote by~$\{e_{\alpha}^*\}$ the basis of
$\ut^*$ dual to the basis $\{e_{\alpha}$, $\alpha\in\Phi^+\}$ of the
algebra $\ut$. Set
\begin{equation*}
f=f_{D,\xi}=\sum_{\beta\in D}\xi_{\beta}e_{\beta}^*\in\ut^*.
\end{equation*}
\defi{We say that the orbit
$\Omega=\Omega_{D,\xi}\subset\ut^*$ of the element $f$ under the
coadjoint action of~the~group $U$ is \emph{associated} with the
subset~$D$. The element~$f$ is called the \emph{canonical form}
on~the~orbit~$\Omega$.} Note that many important examples deal with
orbits associated with {ortho\-gonal} subsets, see
Subsection~\ref{sst_exam}.

\sst To formulate the main result, we need some facts concerning
involutions in the Weyl group of~the~root system $\Phi$. Namely, for
a given orthogonal subset $D\subset\Phi^+$, we put
\begin{equation*}
\sigma=\sigma_D=\prod_{\beta\in D}r_{\beta}\in W
\end{equation*}
(commuting reflections $r_{\beta}$ are taken in any fixed order).
Obviously, $\sigma$ is an involution, i.e., an element of order two
of the group $W$.

To each element $w\in W$ one can assign the numbers $l(w)$ and
$s(w)$. By definition, $l(w)$ (resp.~$s(w)$) is the length of a
reduced (the shortest) expression of $w$ as a product of simple
(resp. arbitrary) reflections. One has $s(\sigma)=|D|$. It's
well-known that $l(\sigma)=|\Phi_{\sigma}|$, where
$\Phi_{\sigma}=\{\alpha\in\Phi^+\mid\sigma\alpha<0\}$. As~usual,
$\alpha>0$ means that $\alpha\in\Phi^+$, and $\alpha<0$ means that
$\alpha\in\Phi^-$. Furthermore, by~$<$ we denote the usual partial
order on $\Phi$: by definition, $\alpha>\beta$ (or $\beta<\alpha$)
if $\alpha-\beta$ is a sum of positive roots.

Things now are ready to formulate the Main Theorem. Since $\Omega$
is an {ir\-re\-du\-cible} affine variety (see \cite[Proposition
8.2]{Humphreys}~and~\cite[Proposition 2.5]{Steinberg2}), one can ask
how to compute $\dim\Omega$, the dimension of $\Omega$ over $k$. (In
fact, if $f$ is an element of $\ut^*(q)$ and $\Omega(q), \Omega$ are
its orbits under the action of~the~groups~$U(q)$, $U$ respectively,
then the complex dimension of the irreducible representation
of~$U(q)$ corresponding to the orbit $\Omega(q)$ is equal to
$q^{\dim\Omega/2}$, see \cite{Kazhdan}.)

\mtheo{Let $D$ be an orthogonal subset of $\Phi^+$\textup{,} $\xi$ a
set of non-zero\linebreak scalars from $k$ and
$\Omega=\Omega_{D,\xi}$ the orbit associated with $D$. Then
$\dim\Omega$ does not depend on $\xi$ and is less or equal
to~$l(\sigma)-s(\sigma)$.\label{mtheo}}

\nota{i) This Theorem proves Conjecture~1.4 from \cite{Ignatev3}.
Note that in many cases (e.g., for~elementary orbits) $\dim\Omega$
is \emph{equal} to $l(\sigma)-s(\sigma)$, see Subsection
\ref{sst_exam}.

ii) On the other hand, for classical groups, the difference between
$\dim\Omega$ and $l(\sigma)-s(\sigma)$ can~be computed explicitly;
furthermore, a polarization of~$\ut$ at the canonical form on
$\Omega$ can be constructed, see~\cite[Theorems 1.1 and
1.2]{Ignatev3} (polarizations play an important role in the explicit
construction of~the~{re\-pre\-sen\-tation} corresponding to a given
orbit). We don't know how to do this for an arbitrary root system.}

\sect{Lemmas and examples} \sst Without loss of generality we can
assume $\Phi$ to be an irreducible root system. Indeed, let
$\Phi=\bigcup_{i=1}^m\Phi_i$ be the decomposition of $\Phi$ into the
union of its pairwise orthogonal irreducible components. Put
$D=\bigcup_{i=1}^mD_i$, where $D_i=D\cap\Phi_i$. Put also
$\ut_i=\sum_{\alpha\in\Phi_i^+}ke_{\alpha}$ for all $i$, and let
$\ut_i^*$ be~the~subspace of~$\ut^*$ dual to the subalgebra $\ut_i$.
Denote by $f_i\in\ut_i^*$ the restriction of~$f$ to~$\ut_i$. Denote
also by $\Omega_i\subset\ut_i^*$ the orbit of $f_i$ under the
coadjoint action of the group $U_i=\exp(\ut_i)$. Finally, denote
by~$W_i$ the~Weyl group of the root system $\Phi_i$ and put
$\sigma_i=\sigma_{D_i}\in W_i$. \lemmp{Suppose that Theorem
$\ref{mtheo}$ holds for all $\Omega_i$, $i=1,\ldots, m$. Then
Theorem $\ref{mtheo}$ holds for the orbit
$\Omega$.\label{lemm_irred}}{Obviously, if~$\alpha\in\Phi_i$,
$\beta\in\Phi_j$ and~$i\neq j$, then~$\alpha+\beta\notin\Phi$. Hence
if~$x=\sum_{i=1}^mx_i$, $y=\sum_{j=1}^my_j$, $x_i,y_i\in\ut_i$,
then~$\ad{-y_j}^rx_i=0$ for all~$j\neq i$, $r>0$. Since
$\exp\ad{-y_i}x_i\in\ut_i$ for all $i$, we obtain
\begin{equation*}
\begin{split}
\exp(y).f(x)&=\left(\sum_{i=1}^mf_i\right)\left(\sum_{j=1}^m\exp
\ad{-y_j}x_j\right)\\
&=\sum_{i=1}^mf_i(\exp\ad{-y_i}x_i)=\sum_{i=1}^m\exp(y_i).f_i(x_i),
\end{split}
\end{equation*}
so the maps
\begin{equation*}
\begin{split}
&\Omega_1\times\ldots\times\Omega_m\to\Omega\colon
(\lambda_1,\ldots,\lambda_m)\mapsto\lambda_1+\ldots+\lambda_m\text{
and}\\
&\Omega\to\Omega_1\times\ldots\times\Omega_m\colon
\lambda\mapsto(\lambda\mathbin{\mid}_{\ut_1},\ldots,\lambda\mathbin{\mid}_{\ut_m})
\end{split}
\end{equation*}
are isomorphisms of affine varieties inverse to each other.

Suppose that Theorem \ref{mtheo} holds for all $\Omega_i$. Let
$\xi_i=(\xi_{\beta})_{\beta\in D_i}$. Since $\Omega_i$ coincides
with $\Omega_{D_i,\xi_i}$, $\dim\Omega$ doesn't depend on $\xi$. On
the other hand, if $i\neq j$, then
$r_{\beta}\mathbin{\mid}_{\Phi_j}=\mathrm{id}_{\Phi_j}$
for~all~$\beta\in D_i$, so ${l(\sigma)=\sum_{i=1}^ml(\sigma_i)}$.
Finally,
$s(\sigma)=|D|=|\bigcup_{i=1}^mD_i|=\sum_{i=1}^m|D_i|=\sum_{i=1}^ms(\sigma_i)$.
This concludes the~proof.}

\newpage
From now on and to the end of the paper, we assume $\Phi$ to be
irreducible.

\sst\label{sst_sing} Sometimes orbits associated with different
orthogonal subsets coincide. To give the precise statement, we need
to introduce the important concept of singular roots.
\defi{Let $\beta\in\Phi^+$ be a positive
root. Roots $\alpha,\gamma\in\Phi^+$ are called
$\beta$-\emph{singular} if $\alpha+\gamma=\beta$. The set of all
$\beta$-singular roots is denoted by~$S(\beta)$.}

Of course, one can easily describe the set of all $\beta$-singular
roots for a given root $\beta$ (see \cite[formula (2)]{Ignatev3} for
the case of classical groups).

Suppose that there exist $\beta_0,\beta_1\in D$ such that
$\beta_0\in S(\beta_1)$. Put $D'=D\setminus\{\beta_0\}$,
${\xi'=\xi\setminus\{\xi_{\beta_0}\}}$, $f'=f_{D',\xi'}$. Let
$\Omega'$ be the orbit of $f'$. \lemmp{The orbit $\Omega$
coincide\footnote {Cf. \cite[Proposition 2.1]{Ignatev3}.} with the
orbit $\Omega'$.\label{lemm_sing}}{Suppose $\beta_1=\beta_0+\alpha$.
Then $\|\alpha\|^2$ is equal to $\|\beta_0\|^2+\|\beta_1\|$. Since
$\Phi$~is~irreducible, we conclude that $\Phi$ is multiply laced
(the root $\alpha$ is long, the roots $\beta_0,\beta_1$ are short);
further, the square of the length of a long root is twice to the
square of the length of a short one. (In~other words, the root
system $\Phi$ is of type $B_n$, $C_n$ or $F_4$.)

Set $\wt f=\exp (ce_{\alpha}).f'$ for some $c\in k^*$. One has
\begin{equation*}
\wt f(e_{\gamma})=f'(e_{\gamma})-c\cdot
f'(\ad{e_{\alpha}}e_{\gamma})+\dfrac{1}{2}c^2\cdot
f'(\ad{e_{\alpha}}^2e_{\gamma})-\ldots
\end{equation*}
for a given root $\gamma\in\Phi^+$. Suppose $\wt
f(e_{\gamma})\neq0$, then there exists $N\gee0$ such that
$\gamma+N\alpha\in D'$. Of course, this holds for $\gamma=\beta_0$
and $N=1$, because $\beta_0+\alpha=\beta_1\in D'$. Suppose that
$N\gee2$ and $\beta_0+N\alpha=\beta\in D'$
(and~so~$\beta\neq\beta_1$). In this case $8||\beta_0||^2\lee
N^2\cdot||\alpha||^2=||\beta_0||^2+||\beta||^2$, a contradiction.
Hence $\wt f(e_{\beta_0})=-c\cdot
N_{\alpha\beta_0}\cdot\xi_{\beta_1}$.

Suppose now that $\gamma\neq\beta_0$ and $\gamma+N\alpha=\beta\in
D'$. If $N=0$, then $\gamma\in D'$, so we can assume ${N\gee1}$. We
see that $\gamma+(N-1)\alpha=\beta-\alpha=\beta-\beta_1+\beta_0$. If
$N=1$, then $\beta\neq\beta_1$,
so~${||\gamma||^2=||\beta+\beta_0||^2+||\beta_1||^2}$. But
$||\beta_1||^2=||\beta_0||^2$ ($\beta_1$~and $\beta_0$ are short),
and $||\beta+\beta_0||^2\gee2||\beta_0||^2$ (the roots
$\beta,\beta_0$ are either equal or~orthogonal). Thus,
$||\gamma||^2\gee3||\beta_0||^2$, a contradiction. Hence~$N\gee2$.
On the other hand,\linebreak
$(\beta,\alpha)=(\beta,\beta_1)-(\beta,\beta_0)$ and
$\gamma=\beta-N\alpha$. If $\beta\neq\beta_1$, then
${(\beta,\alpha)\lee0}$,
so~$||\gamma||^2=||\beta||^2+N^2\cdot||\alpha||^2-2N\cdot(\beta,\alpha)\gee4||\alpha^2||$,
a contradiction. But if ${\beta=\beta_1}$, then $||\gamma||^2=
||N\beta_0+(1-N)\beta_1||^2=(N^2+(1-N)^2)||\beta_0||^2\gee5||\beta_0||^2$,
a contradiction.

We conclude that $\wt f(e_{\gamma})=f'(e_{\gamma})=f(e_{\gamma})$
for all $\gamma\neq\beta_0$. Hence if
$c=-\xi_{\beta_0}/(N_{\alpha\beta_0}\cdot\xi_{\beta_1})$, then $\wt
f$ coincides with $f$. By definition, $\wt f\in\Omega'$, so
$\Omega'=\Omega$ as required.}

From now on and to the end of the paper, we assume that
$S(\beta)\cap D=\varnothing$ for all $\beta\in D$.

\sst\label{sst_zapr} To prove the Main Theorem for simply laced root
systems, we need some more preparations. Let $\eta$, $\eta'$,
$\eta_i$, $\theta$, $\theta'$, $\theta_j$, $\psi$, $\psi'$,
$\psi_l$, $\psi_l'$ be distinct positive roots; assume the roots
$\eta$, $\eta'$, $\eta_i$ to~be~pairwise orthogonal and assume the
root $\eta$ to be maximal among all $\eta$'s w.r.t the usual order
on~$\Phi$. Consider the following cases:
\begin{equation*}
\begin{split}
&\begin{aligned}1.\quad&\eta=\theta+\psi=\theta'+\psi',&
2.\quad&\eta=\theta+\psi,&3.\quad&\eta=\theta+\psi,\\
&\eta_1=\theta_1+\psi,&
&\eta'=\theta'+\psi,& &\eta_1=\theta_1+\psi,\\
&\eta_2=\theta_2+\psi,&
&\eta_1=\theta+\psi_1,& &\eta_2=\theta+\psi_2,\\
&\vphantom{dfrac{1}{\sum}}\eta'=\theta_1+\psi',& &\eta_2=\theta+\psi_2,& &\eta_3=\theta_1+\psi_3,\\
\end{aligned}\\
&\begin{aligned}4.\quad&\eta=\theta+\psi=\theta_1'+\psi_1',\\
&\eta_1=\theta_1+\psi=\theta_1'+\psi',\\
&\eta'=\theta+\psi'=\theta_1+\psi_1',\\
&\{\theta,\theta_1,\theta_1',\psi,\psi',\psi_1'\}\cap
S(\eta_2)\neq\varnothing.\\
\end{aligned}
\end{split}
\end{equation*}

\defi{A set of positive roots
satisfying the conditions of type $1$, $2$, $3$ or $4$ is called
\emph{non-admissible}.}

\lemmp{There are no non-admissible subsets in
$D_5^+$.\label{lemm_zapr}}{Straightforward.}

\sst\label{sst_exam} Before the proof of the Main Theorem, let us
consider some examples of orbits associated with orthogonal subsets.
Let us firstly consider the case $\Phi=A_{n-1}$ (i.e.,
$U=\mathrm{UT}_n$, the unitriangular group). It's convenient to
identify $A_{n-1}^+$ with the subset of $\Rp^n$ of the from
$\{\epsi_i-\epsi_j, 1\lee i<j\lee n\}$ (by~$\{\epsi_i\}_{i=1}^n$ we
denote the standard basis of $\Rp^n$). The Weyl group of $A_{n-1}$
is isomorphic to~$S_n$, the~symmetric group on $n$ letters.

\exam{Let
$D=D_{\reg}=\{\epsi_1-\epsi_n,\epsi_2-\epsi_{n-1},\ldots,\epsi_{n_1}-\epsi_{n_2}\}$,
where ${n_1=[n/2]}$, ${n_2=n-n_1+1}$. Then $\sigma$ is the longest
element of the Weyl group~$W$ and $\Phi_{\sigma}=\Phi^+$ (i.e.,
$\sigma(\alpha)<0$ for~\emph{all} positive roots $\alpha$). Then the
orbit $\Omega$ is \emph{regular}, i.e., has the maximal dimension
among all coadjoint orbits. The dimension of $\Omega$ equals
$$\dim\Omega=l(\sigma)-s(\sigma)=2\mu(n),\text{ }\mu(n)=(n-2)+(n-4)+\ldots$$

ii) Now let
\begin{equation*}
D=D_{\sreg}=(D_{\reg}\setminus\{\epsi_i-\epsi_{n-i+1},
\epsi_{i+1}-\epsi_{n-i}\})\cup\{\epsi_i-\epsi_{n-i},\epsi_{i+1}-\epsi_{n-i+1}\}.
\end{equation*}
In this case, the orbit $\Omega$ is \emph{subregular}, i.e., has the
second maximal dimension $\dim\Omega=l(\sigma)-s(\sigma)=2\mu(n)-2$,
see~\cite[Section 3]{IgnatevPanov}.}

\exam{Let $\Phi$ be an arbitrary root system. Suppose that $|D|=1$.
Then the orbit $\Omega$ is~called \emph{elementary}. It's easy to
see that $\dim\Omega=|S(\beta)|$ \cite[Section 4]{Mukherjee}. It's
straightforward to check that $l(\sigma)-s(\sigma)$ coincides with
$|S(\beta)|$ (see \cite[Section 4]{Ignatev3} for the case of
classical groups).\label{exam_elem}}

\exam{On the other hand, if $\Phi=B_3$ (recall that
$B_3^+=\{\epsi_i$, $\epsi_i\pm\epsi_j$, $1\lee i<j\lee 3\}$) and
$D=\{\epsi_1,\epsi_2+\epsi_3\}$, then the dimension of $\Omega$ is
\emph{less} than $l(\sigma)-s(\sigma)$, because $\dim\Omega=4$ and
$l(\sigma)-s(\sigma)=6$ (see \cite{Ignatev2} or \cite{Ignatev3}).}

\bigskip\sect{Simply laced root systems}

\sst Throughout this Section, $\Phi$ is a \emph{simply laced} root
system, i.e., all roots from $\Phi$ have the same length. (In other
words, $\Phi$ is of type $A_n$, $D_n$, $E_6$, $E_7$ or $E_8$).
Without loss of generality, suppose that the length of a root from
$\Phi$ equals $1$. Then the inner product of two non-orthogonal
roots from $\Phi$ equals either $\pm1$ or $\pm1/2$. Moreover,
suppose $\alpha,\beta\in\Phi^+$, then $(\beta,\alpha)=1/2$ if and
only if either $\alpha\in S(\beta)$ or $\beta\in S(\alpha)$; in this
case, $r_{\beta}\alpha=\alpha-\beta$. On the other hand,
$(\alpha,\beta)=-1/2$ if and only if $\alpha+\beta\in\Phi^+$;
in~this case, $r_{\beta}\alpha=\alpha+\beta$.

As above, let $D$ be an orthogonal subset of $\Phi^+$, $\xi$ a set
of non-zero scalars from $k$, $\Omega=\Omega_{D,\xi}$ the~associated
coadjoint orbit and $f$ the~canonical form on~$\Omega$. Firstly, let
us prove that the dimension of~the orbit $\Omega$ is less or equal
to$l(\sigma)-s(\sigma)$. The proof is by induction on the rank of
$\Phi$. The base ($\rk\Phi=1$, i.e., $\Phi=A_1$) is straightforward.
To perform the inductive step, it's enough to prove the statement
only for irreducible root systems of a given rank, as shows
Lemma~\ref{lemm_irred}.

For the case $|D|=1$ (i.e., the case of elementary orbits), there is
nothing to prove, see Example~\ref{exam_elem}. Suppose $|D|>1$. Pick
a root $\beta$ maximal among all roots from $D$. Put $\wt
D=D\setminus\{\beta\}$. In order to use the inductive hypothesis,
we'll define the root system of rank less than the rank of $\Phi$.
Precisely, put $\Au=\{\alpha\in\Phi^+\mid(\alpha,\beta)\neq0\}$ and
$\wt\Phi=\pm\wt\Phi^+$, where
\begin{equation*}
\wt\Phi^+=\Phi^+\setminus\Au=\{\alpha\in\Phi^+\mid(\alpha,\beta)=0\}.
\end{equation*}
The following Lemma is obvious.\newpage

\lemm{The set $\wt\Phi$ is a root system.\label{post_RS}}

By construction, $\rk\wt\Phi<\rk\Phi$. Obviously, $\wt D=
D\cap\wt\Phi^+$. Denote by~$\wt\ut$ the subalgebra of~$\ut$ spanned
by all vectors~$e_{\alpha}$, $\alpha\in\wt\Phi^+$. Put $\wt
f=f\mathbin\mid_{\wt\ut}\in\wt\ut^*\subset\ut^*$, so
$f=\xi_{\beta}e_{\beta}^*+\wt f$. Put also $\wt U=\exp(\wt\ut)$,
$\wt\xi=\xi\setminus\{\xi_{\beta}\}$. Let $\wt\Omega=\Omega_{\wt
D,\wt\xi}\subset\wt\ut^*$ be the coadjoint orbit of~the~group~$\wt
U$ associated with $\wt D$. Then $\wt f$ is the canonical form
on~the~orbit~$\wt\Omega$.

Finally, let $\wt\sigma$ be the involution in the Weyl group $\wt W$
of the root system $\wt\Phi$ corresponding to~the subset $\wt D$. By
the inductive assumption, $\dim\wt\Omega$ is less or equal to
$l(\wt\sigma)-s(\wt\sigma)$. Obviously, $s(\sigma)=s(\wt\sigma)+1$,
so~it~remains to compare $l(\sigma)$ with $l(\wt\sigma)$.

\sst\label{sst_rad} Let $\at$ be the \emph{radical} of the bilinear
form $x,y\mapsto f([x, y])$, $x, y\in\ut$, i.e.,
\begin{equation*}
\at=\rad_{\ut}f=\{x\in\ut\mid f([x, y])=0\text{ for all }y\in\ut\}.
\end{equation*}
It's well-known that $\dim\Omega=\codim_{\ut}\at=|\Phi^+|-\dim\at$
\cite[Section 3]{AndreNeto}. Similarly, let $\wt\at=\rad_{\wt\ut}\wt
f$ be~the~radical of the bilinear form $x, y\mapsto\wt f([x, y])$,
$x, y\in\wt\ut$. Then
$\dim\wt\Omega=\codim_{\wt\ut}\wt\at=|\wt\Phi^+|-\dim\wt\at$. Put
$\ut_{\Au}=\sum_{\alpha\in\Au}ke_{\alpha}$ (so
$\ut=\ut_{\Au}\oplus\wt\ut$ as vector spaces) and
$\bt=\at\cap\ut_{\Au}$.

\lemmp{The subalgebra $\at$ coincides with the direct sum of its
subspaces $\bt$~and~$\wt\at$, i.e.,
$\at=\bt\oplus\wt\at$.\label{post_oplus}}{i) For any
$x=\sum_{\alpha\in\Phi^+}b_{\alpha}e_{\alpha}\in\ut$ let $\Supp
x=\{\alpha\in\Phi^+\mid b_{\alpha}\neq0\}$. Let
$\rt=\rad_{\ut}\xi_{\beta}e_{\beta}^*=\rad_{\ut}e_{\beta}^*=\langle
e_{\alpha}\mid\alpha\in\Phi^+\setminus
S(\beta)\rangle_k\supset\wt\ut$. Suppose that there exists
$x\in\wt\at$ such that $x\notin\at\cap\wt\ut$. Then there exists
$y\in\ut$ such that $f([x, y])\neq0$. Since $x\in\wt\ut\subset\rt$,
we see that $\xi_{\beta}e_{\beta}^*([x, y])=0$, so $\wt f([x,
y])\neq0$. Hence there exist $\alpha\in\Supp x$,
$\gamma\in\Au\cap\Supp y$ such that $\alpha+\gamma=\wt\beta\in\wt
D$. Using the orthogonality of the subset $D$ and the fact that
$\gamma\in\Au$, we get
$(\alpha,\beta)=(\wt\beta-\gamma,\beta)=(\wt\beta,\beta)-
(\gamma,\beta)=-(\gamma,\beta)\neq0$. This stands in contradiction
with the choice of $\alpha\in\wt\Phi^+$. Thus, $x\in\at\cap\wt\ut$,
so $\wt\at\subset\at\cap\wt\ut\subset\at$.

ii) On the other hand, let $x=y+z\in\at$, where $y\in\ut_{\Au}$,
$z\in\wt\ut$. If $\gamma\in\wt\Phi^+$, then $\gamma\notin S(\beta)$
and $\alpha+\gamma\in\Au$ for all $\alpha\in\Au$. Hence
$f([x,e_{\gamma}])=\xi_{\beta}e_{\beta}^*([x,e_{\gamma}])+\wt
f([y,e_{\gamma}])+\wt f([z,e_{\gamma}])=\wt f([z,e_{\gamma}])=0$,
i.e., $z\in\wt\at$. According to the step i),
$z\in\at\cap\wt\ut\subset\at$. Consequently $y=x-z\in\at$, so
$y\in\at\cap\ut_{\Au}=\bt$ and~$\at=\bt+\wt\at$. But
$\bt\cap\wt\at=0$, so the sum is direct. This concludes the proof.}

\sst To prove the inequality $\dim\Omega\lee l(\sigma)-s(\sigma)$,
we need the following key observation. \lemmp{The inequality
$\#\{\alpha\in\Au\mid\sigma\alpha>0\}+1\lee\dim\bt$
holds.\label{lemm_Au}}{Let
$\wt\Au=\{\alpha\in\Au\mid\sigma\alpha>0\}\cup\{\beta\}$ (clearly,
$\sigma\beta=-\beta<0$). It's enough to construct a linearly
independent set $\{x_{\alpha}\}_{\alpha\in\wt\Au}\subset\bt$. Since
$(\beta,\beta)=1$ and $\beta$ is not singular to any root from $D$,
$\beta\in\Au$ and $x_{\beta}=e_{\beta}\in\bt$.

It's convenient to split the set $\wt\Au$ into a union
$\wt\Au=\Au^+\cup\Au^-\cup\{\beta\}$, where
\begin{equation*}
\begin{split}
&\Au^+=\{\alpha\in\Au\mid\sigma\alpha>0\text{ and
}(\alpha,\beta)>0\},\\
&\Au^-=\{\alpha\in\Au\mid\sigma\alpha>0\text{ and
}(\alpha,\beta)<0\}.
\end{split}
\end{equation*}
Let's consider two different cases, $\alpha\in\Au^-$ and
$\alpha\in\Au^+$.

\medskip i) First, let $\alpha\in\Au^-$, i.e.,
$(\alpha,\beta)=-1/2<0$ and $\sigma\alpha>0$. Suppose $\alpha$~is
singular to the roots $\beta_1,\ldots,\beta_l\in D$ and is not
singular to any other root from~$D$. Put
$\gamma_i=\beta_i-\alpha\in\Phi^+$ for all $i=1,\ldots,l$. Then
\begin{equation*}
(\beta,\gamma_i)=(\beta,\beta_i-\alpha)=-(\beta,\alpha)=1/2,
\end{equation*}
so either $\gamma_i\in S(\beta)$ or $\beta\in S(\gamma_i)$. But if
the second case occurs, then $\beta<\beta_i$, a contradiction with
the choice of the root $\beta$. Thus, $\gamma_i\in S(\beta)$ for all
$i$, Put $\alpha_i=\beta-\gamma_i\in\Phi^+$, $i=1,\ldots, l$.

Now set $b_i=-(\xi_{\beta_i}\cdot
N_{\alpha\gamma_i})/(\xi_{\beta}\cdot N_{\alpha_i\gamma_i})$,
$i=1,\ldots,l$, and
$x_{\alpha}=e_{\alpha}+\sum_{i=1}^lb_ie_{\alpha_i}$. Clearly,
$x_{\alpha}\in\ut_{\Au}$. We claim that
$x_{\alpha}\in\at=\rad_{\ut}f$. Indeed, let $\delta$ be a positive
root. By definition,
\begin{equation*}\predisplaypenalty=0
f([x_{\alpha},e_{\delta}])=N_{\alpha\delta}\cdot
f(e_{\alpha+\delta})+\sum\nolimits_{i=1}^lN_{\alpha_i\delta}\cdot
f(e_{\alpha_i+\delta})\cdot b_i.
\end{equation*}

Pick a number $i$. We note that $\alpha_i+\dl\notin D$ if
$\dl\neq\gamma_i$. Indeed, assume the converse. Then there exists
$\wt\beta\in D$ such that $\wt\beta\neq\beta$ and
$\alpha_i+\dl=\wt\beta$. (Obviously, $\dl\neq\gamma_i$ is equivalent
to $\wt\beta\neq\beta$.) Since $(\wt\beta,\alpha_i)=1/2$,
$(\wt\beta,\gamma_i)=(\wt\beta,\beta-\alpha_i)=-1/2$. Hence
$(\wt\beta,\alpha)=(\wt\beta,\beta_i-\gamma_i)=1/2$. Let
$\beta_1',\ldots,\beta_s'$ be all the roots from $D$, which aren't
orthogonal to $\alpha$ except the roots $\beta,\beta_i,\wt\beta$.
Then
\begin{equation*}
\sigma\alpha=\alpha+\beta-\beta_i+\wt\beta-2(\alpha,\beta_1')
\cdot\beta_1'-\ldots-2(\alpha,\beta_s')\cdot\beta_s'.
\end{equation*}
Since $(\alpha,\beta_r')=\pm1/2$ for all $1\lee r\lee s$, we obtain
\begin{equation*}
||\sigma\alpha-\alpha||^2=2-2(\sigma\alpha,\alpha)=
||\beta||^2+||\beta_i||^2+||\wt\beta||^2+\sum_{r=1}^s||\beta_r'||^2=3+s.
\end{equation*}
We see that either $s=0$ or $s=1$, because
$(\sigma\alpha,\alpha)\gee-1$. If $s=0$, then
\begin{equation*}
\sigma\alpha=\alpha+\beta-\beta_i-\wt\beta=
\alpha+(\alpha_i+\gamma_i)-(\alpha+\gamma_i)-(\alpha_i+\delta)=-\delta<0.
\end{equation*}
On the other hand, $(\alpha,\beta_i)=1/2$, because $\alpha\in
S(\beta_i)$, so if $s=1$, then
\begin{equation*}
\begin{split}
(\sigma\alpha,\alpha)&=||\alpha||^2+(\beta,\alpha)-(\beta_i,\alpha)
-(\wt\beta,\alpha)-2(\beta_1',\alpha)^2\\
&=1-1/2-1/2-1/2-1/2=-1,
\end{split}
\end{equation*}
i.e., ${\sigma\alpha=-\alpha<0}$. By the way, $\sigma\alpha<0$. This
contradicts the choice of $\alpha$. We conclude that
$\alpha_i+\delta\notin D$ if $\delta\neq\gamma_i$.

Hence if $\dl=\gamma_i$ for some $i$, then
\begin{equation*}
\begin{split}
f([x_{\alpha},
e_{\dl}])&=f([x_{\alpha},e_{\gamma_i}])=N_{\alpha\gamma_i}\cdot\xi_{\beta_i}+
N_{\alpha_i\gamma_i}\cdot\xi_{\beta}\cdot
b_i\\
&=N_{\alpha\gamma_i}\cdot\xi_{\beta_i}-
N_{\alpha_i\gamma_i}\cdot\xi_{\beta}\cdot(\xi_{\beta_i}\cdot
N_{\alpha\gamma_i})/(\xi_{\beta}\cdot
N_{\alpha_i\gamma_i})\\
&=N_{\alpha\gamma_i}\cdot\xi_{\beta_i}-N_{\alpha\gamma_i}\cdot\xi_{\beta_i}=0.
\end{split}
\end{equation*}

On the other hand, if $\delta\neq\gamma_i$ for all $i$, then
$\alpha_i+\dl\notin D$, $1\lee i\lee l$, as above. But
$\alpha+\delta\notin D$, so $f([x_{\alpha}, e_{\dl}])=0$ in the
case. Whence for a given $\alpha\in\Au^-$ the vector $x_{\alpha}$
belongs to $\bt=\at\cap\ut_{\Au}$ as required.

\medskip ii) Let us now consider the case $\alpha\in\Au^+$, i.e.,
$(\alpha,\beta)=1/2>0$ and $\sigma\alpha>0$. The Weyl group acts by
orthogonal transformations, so
$(\beta,\sigma\alpha)=(\sigma\beta,\alpha)=(-\beta,\alpha)=-1/2$.
This yields that $\beta+\sigma\alpha\in\Phi^+$. If
$\beta+\sigma\alpha\in S(\wt\beta)$ for a~some root $\wt\beta\in D$,
then $\beta<\wt\beta$. This contradicts the choice of~$\beta$. Thus,
for a given root ${\alpha\in\Au^+}$, the vector
$x_{\alpha}=e_{\beta+\sigma\alpha}$ belongs to $\bt$. Note also that
$(\beta,\beta+\sigma\alpha)=1-1/2=1/2\neq0$, so
$\beta+\sigma\alpha\in\Au$.

\medskip For a given root $\alpha\in\wt\Au$, we constructed the vector
$x_{\alpha}\in\bt$. It remains to check that the vectors
$x_{\alpha}$, $\alpha\in\Au$, are linearly independent. Since
$\beta+\sigma\alpha$, $\alpha\in\Au^+$, are distinct, the
corresponding vectors $x_{\alpha}=e_{\beta+\sigma\alpha}$ are
linearly independent. If $\alpha\in\Au^-$, then
$e_{\alpha}\in\Supp{x_{\alpha}}$ and
$\Supp{x_{\alpha}}\setminus\{\alpha\}\subset S(\beta)$. Consequently
$x_{\alpha}$, $\alpha\in\Au^-$, are linearly independent, too. Their
union with $x_{\alpha}$, $\alpha\in\Au^+$, is also linearly
independent, because  $(\Au^-\cup
S(\beta))\cap(\beta+\sigma\Au^+)=\varnothing$. Indeed, the inner
products of $\beta$ with roots from $\Au^-$ (resp. from
$\beta+\sigma\Au^+$) are negative (resp. positive), so these subsets
are disjoint. Finally, for a given $\alpha\in\Au^+$, the root
$\beta+\sigma\alpha\in\Phi^+$ isn't $\beta$-singular, because
$\beta+\sigma\alpha>\beta$. Thus, the set
$\{x_{\alpha}\}_{\alpha\in\wt\Au}$ is linearly independent. This
completes the proof.}

\sst Now we'll conclude the proof of the inequality $\dim\Omega\lee
l(\sigma)-s(\sigma)$. \propp{Let $\Phi$ be a reduced irreducible
simply laced root system\textup{,} $D\subset\Phi^+$ an orthogonal
subset\textup{,} $\Omega$ an associated orbit of the group
$U$\textup{,} and $\sigma\in W$ the involution corresponding to $D$.
Then $\dim\Omega\lee
l(\sigma)-s(\sigma)$.\label{prop_dim_s_laced}}{By the above (see
Subsection~\ref{sst_rad}) and the inductive hypothesis,
\begin{equation*}
\begin{split}
\dim\Omega&=|\Phi^+|-\dim\at=|\wt\Phi^+|+|\Au|-\dim\at-\dim\wt\at+\dim\wt\at\\
&=(|\wt\Phi^+|-\dim\wt\at)+|\Au|-(\dim\at-\dim\wt\at)\\
&=\dim\wt\Omega+|\Au|-(\dim\at-\dim\wt\at)\\
&\lee
l(\wt\sigma)-s(\wt\sigma)+|\Au|-(\dim\at-\dim\wt\at)=l(\wt\sigma)-s(\wt\sigma)+|\Au|-\dim\bt.
\end{split}
\end{equation*}

It remains to check that $l(\sigma)-s(\sigma)\gee
l(\wt\sigma)-s(\wt\sigma)+|\Au|-\dim\bt$. But
$s(\sigma)=s(\wt\sigma)+1$. We note also that the reflection
$r_{\beta}$ acts on $\wt\Phi^+$ trivially, so
$|\Phi_{\sigma}\cap\wt\Phi^+|=|\wt\Phi_{\wt\sigma}|=l(\wt\sigma)$
and
$$l(\sigma)=|\Phi_{\sigma}|=|\Phi_{\sigma}\cap\wt\Phi^+|+|\Phi_{\sigma}\cap\Au|
=l(\wt\sigma)+\#\{\alpha\in\Au\mid\sigma\alpha<0\}.$$ Hence it's
enough to prove that $$|\Au|-\#\{\alpha\in\Au\mid\sigma\alpha<0\}+1=
\#\{\alpha\in\Au\mid\sigma\alpha>0\}+1\lee\dim\bt,$$ but this
follows immediately from Lemma~\ref{lemm_Au}.}

\sst In the remainder of the Section, we prove that $\dim\Omega$
doesn't depend on~$\xi$. Let $\xi'=(\xi_{\beta}')_{\beta\in D}$ be a
set of non-zero scalars and $f$ the canonical form on the orbit
$\Omega'=\Omega_{D,\xi'}$. Put also $\at'=\rad_{\ut}f'$. Arguing as
in the proof of Proposition~\ref{post_oplus}, we conclude
that~$\at'=\bt'\oplus\wt\at'$ as vector spaces, where\linebreak
$\bt'=\at'\cap\ut_{\Au}$, $\wt f'=f'\mathbin\mid_{\wt\ut}$ and
$\wt\at'=\rad_{\wt\ut}\wt f'$.

\propp{Let $\Phi$ be a reduced irreducible sumply laced root
system\textup{,} $D\subset\Phi^+$ an orthogonal subset\textup{,} and
$\xi,\xi'$ sets of non-zero scalars. Put
$\Omega=\Omega_{D,\xi}$\textup{,} $\Omega'=\Omega_{D,\xi'}$. Then
$\dim\Omega=\dim\Omega'$.\label{prop_xi_s_laced}}{As above,
$\dim\Omega=\codim_{\ut}\at$ and $\dim\Omega'=\codim_{\ut}\at'$, so
it remains to check that $\dim\at=\dim\at'$. We proceed by induction
on the rank of $\Phi$. The base ($\rk\Phi=1$, i.e., $\Phi=A_1$) is
evident. But $\at=\bt\oplus\wt\at$, $\at'=\bt'\oplus\wt\at'$, and
$\dim\wt\at=\dim\wt\at'$ by an inductive assumption, since
$\rk\wt\Phi<\rk\Phi$. Thus, it's enough to show that
$\dim\bt=\dim\bt'$.

Obviously, it's enough to prove that $\dim\bt\lee\dim\bt'$. Set
$x=\sum_{\alpha\in\Au}x_{\alpha}e_{\alpha}\in\bt$. Put
$y=\vfi(x)=\sum_{\alpha\in\Supp{x}}y_{\alpha}e_{\alpha}$. In the
next Subsection we prove that there exist $y_{\alpha}$ such that
$y\in\bt'$ and if the vectors $x_1,\ldots,x_m$ are linearly
independent, then the vectors $\vfi(x_1),\ldots,\vfi(x_m)$ are
linearly independent, too. Applying this to an arbitrary basis $x_i$
of the space $\bt$, we'll obtain the result.}

\sst In this Subsection, we conclude the proof of Proposition
\ref{prop_xi_s_laced}. Our first goal is to determine the
coefficients $y_{\alpha}$. We set $y_{\alpha}=x_{\alpha}$ for all
$\alpha\in\Au$ except the following four cases.

i) There exists $\alpha_0,\gamma\in\Phi^+$, $\beta_0\in D$ such that
\begin{equation*}
\beta=\alpha+\gamma,\quad\beta_0=\alpha_0+\gamma,
\end{equation*}
and $\alpha,\alpha_0$ aren't singular to any other root from $D$.
Then we put $y_{\alpha}=x_{\alpha}\cdot\xi_{\beta}/\xi_{\beta}'$.

ii) There exist $\wt\alpha,\alpha_0,\gamma,\wt\gamma\in\Phi^+$,
$\beta_0,\wt\beta_0\in D$ such that
\begin{equation*}
\beta=\alpha+\gamma=\wt\alpha+\wt\gamma,\quad
\beta_0=\alpha_0+\gamma,\quad\wt\beta_0=\alpha_0+\wt\gamma,
\end{equation*}
and $\alpha,\wt\alpha,\gamma,\wt\gamma$ aren't singular to any other
root from the subset $D$. Here we let
$y_{\alpha}=x_{\alpha}\cdot(\xi_{\beta}\cdot\xi_{\beta_0}')/
(\xi_{\beta}'\cdot\xi_{\beta_0})$. Since the conditions above are
invariant under the interchanging $\alpha$ and $\wt\alpha$, we also
put
$y_{\wt\alpha}=x_{\wt\alpha}\cdot(\xi_{\beta}\cdot\xi_{\wt\beta_0}')/
(\xi_{\beta}'\cdot\xi_{\wt\beta_0})$.

iii) There exist
$\wt\alpha,\alpha_0,\gamma,\wt\gamma,\gamma_0\in\Phi^+$,
$\beta_0,\wt\beta_0\in D$ such that
\begin{equation*}\postdisplaypenalty=0
\begin{split}
&\beta=\alpha+\gamma=\wt\alpha+\wt\gamma,\\
&\beta_0=\alpha_0+\gamma=\wt\alpha+\gamma_0,\\
&\wt\beta_0=\alpha+\gamma_0=\alpha_0+\wt\gamma,\\
\end{split}
\end{equation*}
and $\alpha,\wt\alpha,\alpha_0,\gamma,\wt\gamma,\gamma_0$ aren't
orthogonal to any other root from $D$. As above, we set
$y_{\alpha}=x_{\alpha}\cdot(\xi_{\beta}\cdot\xi_{\beta_0}')/
(\xi_{\beta}'\cdot\xi_{\beta_0})$. Since the conditions are
invariant under the interchanging $\alpha$ and $\wt\alpha$, we also
put
$y_{\wt\alpha}=x_{\wt\alpha}\cdot(\xi_{\beta}\cdot\xi_{\wt\beta_0}')/
(\xi_{\beta}'\cdot\xi_{\wt\beta_0})$.

iv) There exist $\alpha',\gamma'\in\Phi^+$, $\beta_0\in D$ such that
$\alpha=\alpha_0$,
\begin{equation*}
\beta=\alpha'+\gamma',\quad\beta_0=\alpha_0+\gamma',
\end{equation*}
and $\alpha=\alpha_0$, $\alpha'$ aren't singular to any other root
from the subset $D$. Then we let
$y_{\alpha_0}=x_{\alpha_0}\cdot\xi_{\beta_0}/\xi_{\beta_0}'$.

\medskip Let us check that
$y_{\alpha}$ are well-defined. Suppose $\alpha\in\Au$. If
$\alpha\notin\Supp{x}$, then $y_{\alpha}=x_{\alpha}=0$, so let
$\alpha\in\Supp{x}$. If $\alpha$ is not singular to any root
from~$D$, then~$y_{\alpha}=x_{\alpha}$. On the other hand, suppose
$\alpha\in S(\beta)$ (i.e., $\beta=\alpha+\gamma$,
$\gamma\in\Phi^+$). If $\gamma$ is not singular to any other root
from~$D$, then $f([x, e_{\gamma}])=\xi_{\beta}\cdot x_{\alpha}\cdot
N_{\alpha\gamma}\neq0$. This stands in contradiction with the choice
of $x\in\bt\subset\at=\rad_{\ut}f$, so there exist
$\alpha_0\in\Supp{x}$, $\beta_0\in D$ such that $\beta_0\neq\beta$
and $\beta_0=\alpha_0+\gamma$. Suppose there exist $\wt\beta_0\in
D$, $\gamma_0\in\Phi^+$ such that $\wt\beta_0\neq\beta$ and
$\wt\beta_0=\alpha+\gamma_0$. If $\wt\beta_0=\beta_0$, then
\begin{equation*}
(\beta_0,\beta)=(\beta_0,\alpha+\gamma)=(\wt\beta_0,\alpha)+(\beta_0,\gamma)=1/2+1/2=1,
\end{equation*}
a contradiction with the orthogonality of $D$. Hence
$\wt\beta_0\neq\beta_0$.

If $\gamma_0$ isn't singular to any root from the subset $D$ except
$\wt\beta_0$, then $f([x, e_{\gamma_0}])=\xi_{\wt\beta_0}\cdot
x_{\alpha}\cdot N_{\alpha\gamma_0}\neq0$), a contradiction. Whence
there exist $\wh\beta\neq\wt\beta_0$ such that
$\wh\beta=\wh\alpha+\gamma_0$. If $\wh\beta\neq\beta_0$, then
consider the set
$\Psi=\langle\beta,\wh\beta,\beta_0,\wt\beta_0,\alpha\rangle_{\Zp}\cap\Phi$.
Clearly, $\Psi$ is a root system of rank~5, and
$\gamma,\gamma_0,\alpha_0,\wh\alpha\in\Psi$. Further,
$\Psi\cap\Phi^+$ is a~\emph{closed} subset of~$\Psi$, i.e., if
$\zeta_1,\zeta_2\in\Psi\cap\Phi^+$ and $\zeta_1+\zeta_2\in\Psi$,
then $\zeta_1+\zeta_2\in\Psi\cap\Phi^+$. According to \cite[\S16,
Exercise 3]{Humphreys2}, all roots from $\Psi\cap\Phi^+$ are
positive with respect to some subsystem of fundamental roots of
$\Psi$. Thus, without loss of generality, we can assume
$\beta,\wh\beta,\beta_0,\wt\beta_0,\gamma,\gamma_0,\alpha,\alpha_0,\wt\alpha$
belongs to $\Psi^+$.

Since the sum of two roots from different irreducible components of
a root system is not a root, $\Psi$ is an irreducible simply laced
root system. But there are \emph{no} four pairwise orthogonal roots
in $A_5^+$, so $\Psi\cong D_5$. Then the roots $\eta=\beta$,
$\eta_1=\wt\beta_0$, $\eta_2=\beta_0$, $\eta_3=\wh\beta$,
$\theta=\gamma$, $\theta_1=\gamma_0$, $\psi=\alpha$,
$\psi_2=\alpha_0$ and $\psi_3=\wh\alpha$ form a non-admissible
subset of~$D_5^+$ of type 3. This contradicts Lemma \ref{lemm_zapr}.

We conclude that
$\wh\beta=\beta_0=\alpha_0+\gamma=\wh\alpha+\gamma_0$. Denote
$\wt\alpha=\wh\alpha$, so
$\beta_0=\alpha_0+\gamma=\wt\alpha+\gamma_0$. Then
$\wt\alpha=\beta_0-\gamma_0=\beta_0-(\wt\beta_0-\alpha)=\beta_0-\wt\beta_0+\alpha$,
so $(\beta,\wt\alpha)=(\beta,\alpha)=1/2$. Hence $\wt\alpha\in
S(\beta)$, because if $\beta\in S(\wt\alpha)$, then $\beta<\beta_0$,
and $\beta$ is not maximal among all roots from $D$. In other words,
$\beta=\wt\alpha+\wt\gamma$ for some~${\wt\gamma\in\Phi^+}$. We see
that
\begin{equation*}
\wt\beta_0=\alpha+\gamma_0=(\beta-\gamma)+(\beta_0-\wt\alpha)=
(\beta_0-\gamma)+(\beta-\wt\alpha)=\alpha_0+\wt\gamma.
\end{equation*}
Lemma \ref{lemm_zapr} shows that
$\alpha,\wt\alpha,\alpha_0,\gamma,\wt\gamma,\gamma_0$ aren't
orthogonal to any other root $\beta_2\in D$. Indeed, assume the
converse, Then the roots $\eta=\beta$, $\eta_1=\beta_0$,
$\eta'=\wt\beta_0$, $\theta=\alpha$, $\theta_1=\alpha_0$,
$\theta_1'=\wt\alpha$, $\psi=\gamma$, $\psi'=\gamma_0$,
$\psi_1'=\wt\gamma$ and~$\eta_2=\beta_2$ form a non-admissible
subset of~$D_5^+$ of type~4. Thus, $\alpha$ belongs to case~iii),
and the roots
$\wt\alpha,\alpha_0,\beta_0,\wt\beta_0,\gamma,\wt\gamma,\gamma_0$
are determined uniquely.

Suppose now that $\beta=\alpha+\gamma$, $\beta_0=\alpha_0+\gamma$,
but~$\alpha$ is not singular to any other root from~$D$ except
$\beta$. Suppose also that there exists $\wt\gamma\in\Phi^+$ such
that $\wt\beta_0=\alpha_0+\wt\gamma$. Then
$$(\beta,\wt\gamma)=(\beta,\wt\beta_0-\alpha_0)=(\beta,\wt\beta_0-\beta_0+\gamma)=1/2,$$
because $\gamma\in S(\beta)$. This implies $\wt\gamma\in S(\beta)$,
because if $\beta\in S(\gamma)$, then~$\beta<\beta_0$, and $\beta$
is not maximal. Let $\beta=\wt\alpha+\wt\gamma$ for some
$\wt\alpha\in\Phi^+$. We note that $\wt\alpha$ is not singular to
any other root from $D$ except $\beta$. Indeed, if the converse
holds, then the root~$\wt\alpha$ belongs to case~iii). But this
yields $\alpha\in S(\wt\beta_0)$, a contradiction.

Besides, $\gamma$ isn't singular to any other root from $D$ except
$\beta$ and $\beta_0$. Indeed, if there exist $\beta_2\in D$,
$\alpha_2\in\Phi^+$ such that $\beta_2\neq\beta$,
$\beta_2\neq\beta_0$ and $\beta_2=\alpha_2+\gamma$, then the roots
$\eta=\beta$, $\eta_1=\beta_0$, $\eta'=\wt\beta_0$,
$\eta_2=\beta_2$, ${\theta=\alpha}$, $\theta'=\wt\alpha$,
$\theta_1=\alpha_0$, $\theta_2=\alpha_2$, $\psi=\gamma$
and~$\psi'=\wt\gamma$ form a non-admissible subset of~$D_5^+$ of
type~1. This stands in contradiction with Lemma~\ref{lemm_zapr}.
Similarly, $\wt\gamma$ isn't singular to any other root from~$D$
except $\beta$,~$\wt\beta_0$. We see that~$\alpha$ belongs to case
ii), and the roots $\beta_0$,~$\wt\beta_0$ are determined uniquely.

Suppose now that $\beta=\alpha+\gamma$, $\beta_0=\alpha_0+\gamma$
and the roots $\alpha$,~$\alpha_0$ aren't singular to any other root
from~$D$ except $\beta$, $\beta_0$ respectively, Then $\alpha$
belongs to case i), and the root $\beta_0$ is determined uniquely.

It remains to consider the case when $\alpha=\alpha_0$ isn't
singular to $\beta$, but there exist $\beta_0\in D$,
$\gamma'\in\Phi^+$ such that $\beta_0\neq\beta$ and
$\beta_0=\alpha_0+\gamma'$. Since $\alpha\notin S(\beta)$,
$\alpha$~doesn't belong to cases i)--iii). If $\gamma'\notin
S(\beta)$, then $\alpha=\alpha_0$ doesn't belong to case~iv), too,
so $y_{\alpha_0}=x_{\alpha_0}$.

Suppose now $\gamma'\in S(\beta)$, i.e., there exists
$\alpha'\in\Phi^+$ such that $\beta=\alpha'+\gamma'$. If
$\alpha'\in\Supp{x}$, then the root $\alpha'$ belongs to one of
cases~i)--iii). This implies $y_{\alpha_0}=x_{\alpha_0}$. If
$\alpha'\notin\Supp{x}$, then there exist $\beta_1\in D$,
$\alpha_1\in\Phi^+$ such that $\beta\neq\beta$, $\beta_1\neq\beta_0$
and $\beta_1=\alpha_1+\gamma'$  (if the converse holds, then $f([x,
e_{\gamma'}])=\xi_{\beta_0}\cdot x_{\alpha_0}\cdot
N_{\alpha_0\gamma'}\neq0$, a contradiction). Hence $\alpha'$,
$\alpha=\alpha_0$ aren't singular to any other root from~$D$ except
$\beta, \beta_0$ respectively.

Indeed, suppose there exist $\wt\beta_0\in D$, $\wt\gamma\in\Phi^+$
such that $\wt\beta_0\neq\beta_0$ and
$\wt\beta_0=\alpha_0+\wt\gamma$. The root $\wt\beta_0$ doesn't
coincide with $\beta$, because $\alpha=\alpha_0$ isn't singular to
$\beta$. If $\wt\beta_0$ coincides with $\beta_1$, then
\begin{equation*}
(\beta_0,\beta_1)=(\alpha_0,\gamma',\beta_1)=(\alpha_0,\wt\beta_0)+(\gamma',\beta_1)
=1/2+1/2=1.
\end{equation*}
At the same time
$\wt\gamma=\wt\beta_0-\alpha_0=\wt\beta_0-\beta_0+\gamma'$, so
$(\beta,\wt\gamma)=1/2$ and $\wt\gamma\in S(\beta)$ (if $\beta\in
S(\wt\gamma)$, then $\beta<\beta_0$, so $\beta$ isn't maximal).
However if $\beta=\wt\alpha+\wt\gamma$, $\wt\alpha\in\Phi^+$, then
the roots $\eta=\beta$, $\eta_1=\beta_0$, $\eta_2=\beta_1$,
$\eta'=\wt\beta_0$, $\theta=\alpha'$, $\theta_1=\alpha=\alpha_0$,
$\theta_2=\alpha_1$, $\psi=\gamma'$ and~$\psi'=\wt\gamma$ form a
non-admissible subset of~$D_5^+$ of type~1, a~contradiction with
Lemma~\ref{lemm_zapr}. This contradiction shows that
$\wt\beta_0\neq\beta_1$.

On the other hand, suppose that there exist $\wt\beta\in D$,
$\wt\gamma\in\Phi^+$ such that $\wt\beta\neq\beta$ and
$\wt\beta=\alpha'+\wt\gamma$. If $\wt\beta$ coincides with
$\beta_0$, then $(\beta,\beta_0)=1$, so $\wt\beta\neq\beta_0$; for
the same reason, $\wt\beta\neq\beta_1$. It follows that the roots
$\eta=\beta$, $\eta'=\wt\beta$, $\eta_1=\beta_0$, $\eta_2=\beta_1$,
$\theta=\gamma$, $\theta'=\wt\gamma$, $\psi=\alpha'$,
$\psi_1=\alpha=\alpha_0$ and~$\psi_2=\alpha_1$ form a non-admissible
subset of~$D_5^+$ of type~2. This contradicts Lemma~\ref{lemm_zapr}.

We have proved that if $\alpha=\alpha_0\notin S(\beta)$, $\alpha\in
S(\beta_0)$ for some $\beta_0\in D$, $\beta_0\neq\beta$, and
$y_{\alpha}\neq x_{\alpha}$, then $\alpha=\alpha_0$ belongs to
case~iv); in particular $\beta_0$ is determined uniquely.

Therefore if $y_{\alpha}\neq x_{\alpha}$, then $\alpha$ belongs to
one of cases i)--iv), and the root $\beta_0$ is determined uniquely.
Since $y_{\alpha}$ depend only on $\beta_0$, they are well-defined.
Denote by $X$, $Y$ the ($|\Au|\times l$)-matrices whose columns
consist of the coordinates of the vectors $x_1,\ldots,x_l\in\bt$
and~$y_1=\vfi(x_1),\ldots$, $y_l=\vfi(x_l)$ respectively in the
basis $\{e_{\alpha}\}_{\alpha\in\Au}$ of the algebra $\ut_{\Au}$.
Let~$T$ be the diagonal ($|\Au|\times|\Au|$)-matrix whose $(\alpha,
\alpha)$-th element equals
\begin{equation*}\predisplaypenalty=0
t_{\alpha, \alpha}=\begin{cases}\xi_{\beta}/\xi_{\beta}',&\text{if
$\alpha$ belongs to case i),}\\
(\xi_{\beta}\cdot\xi_{\beta_0}')/
(\xi_{\beta}'\cdot\xi_{\beta_0}),&\text{if
$\alpha$ belongs either to cases ii) or iii),}\\
\xi_{\beta_0}/\xi_{\beta_0}',&\text{if
$\alpha=\alpha_0$ belongs to case iv),}\\
1&\text{otherwise.}\\
\end{cases}
\end{equation*}
We see that $Y=TX$, but $\det T\neq0$, so $\rk{X}=\rk{Y}$. To
conclude the proof, it remains to check that if $x\in\bt$, then
$y\in\bt'$, i.e., $f'([y, e_{\wh\gamma}])=0$ for all
$\wh\gamma\in\Phi^+$. Let us consider four cases.

\medskip 1. Firstly, suppose $(\beta,\wh\gamma)=0$, i.e.,
$\wh\gamma\in\wt\Phi^+$ and $\wh\gamma\notin S(\beta)$. Let
$\wh\gamma$ be singular to the roots $\beta_1,\ldots,\beta_l$
from~$\wt D=D\setminus\{\beta\}$ and not singular to any other root
from $\wt D$. Denote $\alpha_i=\beta_i-\wh\gamma$. Then
$(\beta,\alpha_i)=(\beta,\beta_i-\wh\gamma)=0$, so $e_{\alpha_i}$
doesn't belong to $\Supp{x}\subset\Au=\Phi^+\setminus\wt\Phi^+$.
Thus, $f'([y,e_{\wh\gamma}])=0$.

\medskip 2. Secondly, suppose $(\beta,\wh\gamma)\neq0$, i.e.,
$\wh\gamma\in\Au$. If $(\beta,\wh\gamma)=1$, then $\wh\gamma=\beta$.
But $\beta$ isn't singular to any root from $D$, so $f'([y,
e_{\wh\gamma}])=f'([y, e_{\beta}])=0$. If $(\beta,\wh\gamma)=-1/2$,
then $\wh\gamma$ isn't singular to $\beta$; in this case, denote
$\gamma_0=\wh\gamma$ (so we must prove that $f'([y,
e_{\gamma_0}])=0$). If $\alpha=\wt\beta_0-\gamma_0\notin \Supp{x}$
for all ${\wt\beta_0\in\wt D}$ such that $\gamma_0\in S(\beta_0)$,
then ${f'([y, e_{\gamma_0}])=0}$. On the other hand, suppose there
exist $\wt\beta_0\in\wt D$, $\alpha\in\Supp{x}$ such that
$\wt\beta_0=\alpha+\gamma_0$.

Since $(\beta, \alpha)=(\beta,\wt\beta_0-\gamma)=1/2$, $\alpha\in
S(\beta)$ (if $\beta\in S(\alpha)$, then $\beta<\wt\beta_0$, so
$\beta$~isn't maximal). In other words, there exists
$\gamma\in\Phi^+$ such that $\beta=\alpha+\gamma$. If $\gamma$ isn't
singular to any other root from $D$, then $f([x,
e_{\gamma}])=\xi_{\beta}\cdot x_{\alpha}\cdot
N_{\alpha\gamma}\neq0$, a contradiction. Whence there exist
$\beta_0\in D$, $\beta_0\neq\beta$, $\alpha_0\in\Phi^+$ such that
$\beta_0=\alpha_0+\gamma$. Arguing as above, we see that $\alpha$
belongs to case~iii), and, consequently,
\begin{equation*}\predisplaypenalty=0
\begin{split}
f'([y,e_{\gamma_0}])&=\xi_{\beta_0}'\cdot N_{\wt\alpha\gamma_0}\cdot
y_{\wt\alpha}+\xi_{\wt\beta_0}'\cdot
N_{\alpha\gamma_0}\cdot y_{\alpha}\\
&=\xi_{\beta_0}'\cdot N_{\wt\alpha\gamma_0}\cdot
x_{\wt\alpha}\cdot(\xi_{\beta}\cdot\xi_{\wt\beta_0}')/
(\xi_{\beta}'\cdot\xi_{\wt\beta_0})\\
&+\xi_{\wt\beta_0}'\cdot N_{\alpha\gamma_0}\cdot
y_{\alpha}\cdot(\xi_{\beta}\cdot\xi_{\beta_0}')/
(\xi_{\beta}'\cdot\xi_{\beta_0})\\
&=(\xi_{\beta_0}\cdot N_{\wt\alpha\gamma_0}\cdot
x_{\wt\alpha}+\xi_{\wt\beta_0}\cdot N_{\alpha\gamma_0}\cdot
x_{\alpha})\cdot(\xi_{\beta_0}'\cdot\xi_{\wt\beta_0}'\cdot\xi_{\beta})/
(\xi_{\beta_0}\cdot\xi_{\wt\beta_0}\cdot\xi_{\beta}')\\
&=f([x,
e_{\gamma_0}])\cdot(\xi_{\beta_0}'\cdot\xi_{\wt\beta_0}'\cdot\xi_{\beta})/
(\xi_{\beta_0}\cdot\xi_{\wt\beta_0}\cdot\xi_{\beta}')=0.
\end{split}
\end{equation*}

3. Thirdly, suppose $(\beta,\wh\gamma)=1/2$ and $\wh\gamma\in
S(\beta)$; in this case, denote $\gamma=\wh\gamma$ (so we must prove
that $f'([y, e_{\gamma}])=0$). Let $\beta=\alpha+\gamma$,
$\alpha\in\Phi^+$. Let also~$\gamma$~be singular to the roots
$\beta_1,\ldots,\beta_l$ from~$\wt D$ and not singular to any other
root from~$\wt D$. Denote $\alpha_i=\beta_i-\gamma$. If
$\alpha\in\Supp{x}$, then, arguing as above, we conclude that
$\alpha$ belongs to one of cases i)--iii). If $\alpha$ belongs to
case i), then all the roots $\alpha_i$ belong to case iv). Hence
\begin{equation*}\predisplaypenalty=0
\begin{split}
f'([y, e_{\gamma}])&=\xi_{\beta}'\cdot N_{\alpha\gamma}\cdot
y_{\alpha}+\sum\nolimits_{i=1}^l\xi_{\beta_i}'\cdot N_{\alpha_i\gamma}\cdot y_{\alpha_i}\\
&=\xi_{\beta}'\cdot N_{\alpha\gamma}\cdot
x_{\alpha}\cdot\xi_{\beta}/\xi_{\beta}'+
\sum\nolimits_{i=1}^l\xi_{\beta_i}'\cdot N_{\alpha_i\gamma}\cdot
x_{\alpha_i}
\cdot\xi_{\beta_i}/\xi_{\beta_i}'\\
&=\xi_{\beta}\cdot N_{\alpha\gamma}\cdot
x_{\alpha}+\sum\nolimits_{i=1}^l\xi_{\beta_i}\cdot
N_{\alpha_i\gamma}\cdot x_{\alpha_i}
=f([x, e_{\gamma}])=0.\\
\end{split}
\end{equation*}

On the other hand, if $\alpha$ belongs either to case ii) or iii),
then $\gamma$ is not singular to any other root from~$D$ except
$\beta$ and $\beta_0=\alpha_0+\gamma$, $y_{\alpha_0}=x_{\alpha_0}$,
and
\begin{equation*}\predisplaypenalty0
\begin{split}
f'([y, e_{\gamma}])&=\xi_{\beta}'\cdot N_{\alpha\gamma}\cdot
y_{\alpha}+\xi_{\beta_0}'\cdot N_{\alpha_0\gamma}\cdot
y_{\alpha_0}\\
&=\xi_{\beta}'\cdot N_{\alpha\gamma}\cdot
x_{\alpha}\cdot(\xi_{\beta}\cdot\xi_{\beta_0}')/
(\xi_{\beta'}\cdot\xi_{\beta_0})+\xi_{\beta_0}'\cdot
N_{\alpha_0\gamma}\cdot
x_{\alpha_0}\\
&=(\xi_{\beta}\cdot N_{\alpha\gamma}\cdot
x_{\alpha}+\xi_{\beta_0}\cdot N_{\alpha_0\gamma}\cdot
x_{\alpha_0})\cdot\xi_{\beta_0}'/\xi_{\beta_0}\\
&=f([x, e_{\gamma}])\cdot\xi_{\beta_0}'/\xi_{\beta_0}=0.\\
\end{split}
\end{equation*}

Assume now that $\alpha\notin\Supp{x}$. If $\alpha_i\notin\Supp{x}$
for all $i$, then $f'([y, e_{\gamma}])=0$, because
$\Supp{x}=\Supp{y}$. At the same time if $\alpha_i\in\Supp{x}$, then
there exists $\alpha_j$ such that $i\neq j$ and
$\alpha_j\in\Supp{x}$. We claim that the root $\alpha_i$ isn't
singular to any other root from $D$ except $\beta_i$. Indeed, assume
$\wt\beta_i=\alpha_i+\wt\gamma$ for some $\wt\beta_i\in D$,
$\wt\gamma\in\Phi^+$. Since
$(\beta,\alpha_i)=(\beta,\beta_i-\gamma)=-1/2$, $\alpha_i\notin
S(\beta)$, so $\wt\beta_i\neq\beta$. It follows that
$(\beta,\wt\gamma)=(\beta,\wt\beta_i-\alpha_i)=1/2$, so
$\wt\gamma\in S(\beta)$ (if $\beta\in S(\wt\gamma)$, then
$\beta<\wt\beta_i$, so $\beta$ isn't maximal). Put
$\beta=\wt\alpha+\wt\gamma$, $\wt\alpha\in\Phi^+$. If
$\wt\beta_i=\beta_j$, then
\begin{equation*}
(\beta_i,\beta_j)=(\alpha_i+\gamma,\beta_j)=(\alpha_i,\wt\beta_i)+(\gamma,
\beta_j)=1/2+1/2=1.
\end{equation*}
Hence $\wt\beta_i\neq\beta_j$, so the roots $\eta=\beta$,
$\eta_1=\beta_i$, $\eta_2=\beta_j$, $\eta'=\wt\beta_i$,
$\theta=\alpha$, $\theta'=\wt\alpha$, $\theta_1=\alpha_i$,
$\theta_2=\alpha_j$, $\psi=\gamma$ and~$\psi'=\wt\gamma$ form a
non-admissible subset of~$D_5^+$ of type~1, a contradiction with
Lemma~\ref{lemm_zapr}. This contradiction shows that $\alpha_i$
isn't singular to any other root from $D$ except $\beta_i$, as
required.

Similarly, suppose there exist $\wt\beta\in D$, $\wt\gamma\in\Phi^+$
such that $\wt\beta\neq\beta$ and $\wt\beta=\alpha+\wt\gamma$. If
$\wt\beta$ coincides with $\beta_i$, then
\begin{equation*}
(\beta,\beta_i)=(\alpha+\gamma,\beta_i)=
(\alpha,\wt\beta)+(\gamma,\beta_i)=1/2+1/2=1.
\end{equation*}
Hence $\wt\beta\neq\beta_i$. For the same reason,
$\wt\beta\neq\beta_j$. But this yields that the roots $\eta=\beta$,
$\eta'=\wt\beta$, $\eta_1=\beta_i$, $\eta_2=\beta_j$,
$\theta=\gamma$, $\theta'=\wt\gamma$, $\psi=\alpha$,
$\psi_1=\alpha_i$ and~$\psi_2=\alpha_j$ form a non-admissible subset
of~$D_5^+$ of type~2. This contradicts Lemma~\ref{lemm_zapr}. We've
proved that $\alpha$ and all~$\alpha_i\in\Supp{x}$ aren't singular
to any other root from~$D$ except $\beta$, $\beta_i$ respectively.
This implies that all $\alpha_i\in\Supp{x}$ belong to case~iv), so
\begin{equation*}\predisplaypenalty0
\begin{split}
f'([y, e_{\gamma}])&=\sum\nolimits_{i=1}^l\xi_{\beta_i}'\cdot
N_{\alpha_i\gamma}\cdot
y_{\alpha_i}=\sum\nolimits_{i=1}^l\xi_{\beta_i}'\cdot
N_{\alpha_i\gamma}\cdot x_{\alpha_i}\cdot\xi_{\beta_i}/\xi_{\beta_i}'\\
&=\sum\nolimits_{i=1}^l\xi_{\beta_i}\cdot N_{\alpha_i\gamma}\cdot
x_{\alpha_i}=f([x, e_{\gamma}])=0.
\end{split}
\end{equation*}

4. Finally, suppose $(\beta,\wh\gamma)=1/2$ and $\beta\in
S(\wh\gamma)$. Then $\wh\gamma$ isn't singular to any root from~$D$,
so $f'([y, e_{\wh\gamma}])=0$. This concludes the proof of
Proposition~\ref{prop_xi_s_laced}.

\nota{The case $\Phi=A_n$ was considered by A.N. Panov in the
paper~\cite{Panov}. The case $\Phi=D_n$ was considered by the author
in the paper~\cite{Ignatev3}. Actually, the new result is obtained
only for the root systems of types $E_6,E_7,E_8$. However note that
the \emph{proofs} are similar for all simply laced root systems.}

\sect{Multiply laced root systems}

\sst Throughout the section, we assume $\Phi$ to be reduced
irreducible \emph{multiply laced} root system (i.e., containing long
and short roots). The cases of $B_n$ and $C_n$ were considered by
the author in the paper~\cite{Ignatev3}, so we'll assume that $\Phi$
is of type $F_4$ or $G_2$. Firstly, suppose that $\Phi=G_2$  (this
case is quite easy).

Recall that
$G_2^+=\{\alpha_1,\alpha_2,\alpha_1+\alpha_2,2\alpha_1+\alpha_2,
2\alpha_1+\alpha_2,3\alpha_1+2\alpha_2\}$, where $||\alpha_1||^2=1$,
$||\alpha_2||^2=3$ and the angle between the vectors
$\alpha_1,\alpha_2$ equals $5\pi/6$. Let $D$ be an orthogonal subset
of $\Phi^+$. Of course, $|D|\leq2$. For $|D|=1$, there is nothing to
prove, because $\Omega$ is an elementary orbit
(see~Example~\ref{exam_elem}). There are three orthogonal subsets of
$G_2^+$ of cardinality two; we'll consider all of them subsequently.
Note that $l(\sigma)=6$  and $l(\sigma)-s(\sigma)=6-2=4$, because
$\sigma$ is the central symmetry. Denote $D=\{\beta_1,\beta_2\}$.

i) $\beta_1=\alpha_1,\beta_2=3\alpha_1+2\alpha_2$. The root
$\alpha_1$ is fundamental, so $S(\alpha_1)=\varnothing$. At the same
time $\beta_2$-singular roots are the following:
$S(\beta_2)=\{\alpha_1,3\alpha_1+\alpha_2\}\cup\{\alpha_1+
\alpha_2,2\alpha_1+\alpha_2\}$. Put
$\Mo=\{\alpha_2,\alpha_1+\alpha_2\}$, $\Po=\Phi^+\setminus\Mo$ and
$\pt=\sum_{\alpha\in\Po}ke_{\alpha}$. One can see that
$\pt\subset\ut$ is an \emph{isotropic} subspace, i.e., $f([x, y])=0$
for all $x,y\in\pt$. (Indeed, if $1\leq i\leq2$ and $\alpha,\gamma$
are $\beta_i$-singular, then $\Po$ doesn't contain both of them.)
Further, if $x\notin\pt$, then $\Supp{x}$ contains at least one of
the roots $\gamma_1=\alpha_2,\gamma_2=\alpha_1+\alpha_2$. Actually
if $x=x_1e_{\gamma_1}+\ldots$, $x_1\neq0$, then $f([x,
e_{3\alpha_1+\alpha_2}])=\xi_{\beta_2}\cdot x_1\cdot
N_{\gamma_1,3\alpha_1+\alpha_2}\neq0$. Similarly, if
$x=x_2e_{\gamma_2}+\ldots$, $x_2\neq0$, then $f([x,
e_{2\alpha_1+\alpha_2}])=\xi_{\beta_2}\cdot x_2\cdot
N_{\gamma_2,2\alpha_1+\alpha_2}\neq0$. Thus, $\pt$ is a
\emph{maximal} isotropic subspace with respect to~the~inclusion
order. Hence $\dim\Omega$ doesn't depend on~$\xi$ and equals
$2\cdot\codim_{\ut}\pt=4=l(\sigma)-s(\sigma)$ (see,
f.e.,~\cite[Section 3]{AndreNeto}).

ii) $\beta_1=\alpha_1+\alpha_2,\beta_2=3\alpha_1+\alpha_2$. Here
$S(\beta_1)=\{\alpha_1,\alpha_2\}$, $S(\beta_2)=\{\alpha_1$,
$2\alpha_1+\alpha_2\}$. Put $\Mo=\{\alpha_1\}$. Let $\Po, \pt$ be
defined as above. Evidently, $\pt$ is an isotropic subspace. On the
other hand, if $x=x_1e_{\alpha_1}+\ldots\notin\pt$, then
$f([x,e_{\alpha_2}])=\xi_{\beta_1}\cdot x_1\cdot
N_{\alpha_1\alpha_2}\neq0$, so $\pt$ is a maximal isotropic subspace
and, consequently, $\dim\Omega$ doesn't depend on $\xi$ and equals
$2\cdot\codim_{\ut}\pt=2<4=l(\sigma)-s(\sigma)$.

iii) $\beta_1=\alpha_2,\beta_2=2\alpha_1+\alpha_2$. The root
$\alpha_2$ is fundamental, so $S(\alpha_2)=\varnothing$. At the same
time $S(\beta_2)=\{\alpha_1,\alpha_1+\alpha_2\}$. Putting
$\Mo=\{\alpha_1\}$, we see that $\dim\Omega$ doesn't depend on $\xi$
and equals $2<4=l(\sigma)-s(\sigma)$, as in the previous step.

\sst Let us now consider the more complicated case $\Phi=F_4$.
Recall that
\begin{equation*}F_4^+=\{\epsi_i,\epsi_i\pm\epsi_j,
(\epsi_1\pm\epsi_2\pm\epsi_3\pm\epsi_4)/2,1\lee
i<j\lee4\}\quad(\text{signs are independent}),
\end{equation*}
where $\alpha_1=\epsi_2-\epsi_3$, $\alpha_2=\epsi_3-\epsi_4$,
$\alpha_3=\epsi_4$, $\alpha_4=(\epsi_1-\epsi_2-\epsi_3-\epsi_4)/2$
are fundamental roots (here $\{\epsi_i\}_{i=1}^4$ is the standard
basis of $\Rp^4$). For convenience, put
$\wt\Phi^+=\{\epsi_i,\epsi_i\pm\epsi_j,1\lee i<j\lee4\}$ and
$\Bu=\Phi^+\setminus\wt\Phi^+$. One has $\wt\Phi\cong B_4$ as root
systems, where $\wt\Phi=\pm\wt\Phi^+$.

We begin with the case $D\subset\wt\Phi^+$. We denote by~$\wt W$ the
Weyl group of the root system $\wt\Phi$. We also denote by
$\wt\sigma$ the involution in the $\wt W$ corresponding to the
subset $D$. Clearly, $s(\wt\sigma)=s(\sigma)=|D|$ and
$l(\wt\sigma)\lee l(\sigma)$. Precisely,
$\Fo=\wt\Fo+\#\{\alpha\in\Bu\mid\sigma\alpha<0\}$ (here we put
$\Fo=l(\sigma)-s(\sigma)$ and $\wt\Fo=l(\wt\sigma)-s(\wt\sigma)$).

As above, denote $\wt\ut=\sum_{\alpha\in\wt\Phi^+}ke_{\alpha}$,
$\ut_{\Bu}=\sum_{\alpha\in\Bu}ke_{\alpha}$ (hence
$\ut=\wt\ut\oplus\ut_{\Bu}$ as vector spaces) and set $\wt
f=f\mathbin\mid_{\wt\ut}$, $\wt U=\exp(\wt\ut)$. Let
$\wt\Omega\subset\wt\ut^*$ be the orbit of $\wt f$ under the
coadjoint action of the group $\wt U$. Let $\at=\rad_{\ut}f$,
$\wt\at=\rad_{\wt\ut}\wt f$. It follows from \cite[Theorem
1.2]{Ignatev3} that $\dim\wt\Omega=\wt\Fo-\teta$, where $\teta$
depends only on $D$, not on $\xi$. Finally, put
$\bt=\at\cap\ut_{\Bu}$.

\lemmp{One has $\at=\wt\at\oplus\bt$ as vector spaces \textup{(}cf.
Lemma $\ref{post_oplus}$\textup{)}.}{Suppose $x\in\wt\at$,
$\alpha\in\Supp{x}$. If $\gamma\in\wt\Phi$, then
$f([x,e_{\gamma}])=0$, because $\wt\at=\rad_{\wt\ut}f$ is the
radical of $f$. If $\gamma\in\Bu$, then $\alpha+\gamma\in\Bu$,
because $f([e_{\alpha},e_{\gamma}])=0$. We see that
$f([x,e_{\gamma}])=0$ for all $\gamma\in\Phi^+$, hence $x\in\at$.
Thus, $\wt\at\subset\at$. On the other hand, suppose $x=y+z\in\at$,
$y\in\wt\ut$, $z\in\ut_{\Bu}$, $\alpha\in\Supp{z}$ and
$\gamma\in\wt\Phi^+$. Then $\alpha+\gamma\in\Bu$, so
$f([z,e_{\gamma}])=0$ and $f([y,e_{\gamma}])=0$, i.e.,
$y\in\wt\at\subset\at$. Hence $z\in\at$ and $\at=\wt\at+\bt$. But
$\wt\at\cap\bt=0$. The proof is complete.}

\lemmp{One has $\bt=\langle
e_{\alpha},\alpha\in\Bu\mid\sigma\alpha>0\rangle_k$ \textup{(}cf.
Lemma~$\ref{lemm_Au}$\textup{)}.}{Set
$\wt\Bu=\{\alpha\in\Bu\mid\sigma\alpha>0\}$. Firstly, suppose that
$D$ doesn't contain the roots $\epsi_1$, $\epsi_1\pm\epsi_j$,
$j=2,3,4$. Then $\sigma\alpha=\epsi_1/2\pm\ldots>0$ for all
$\alpha\in\Bu$, so $\wt\Bu=\Bu$. In this case,
$\Bu\cap\bigcup_{\beta'\in D}S(\beta')=\varnothing$. Hence
$\bt=\ut_{\Bu}$ as required.

Secondly, suppose $\beta=\epsi_1\in D$. Then $D$ doesn't contain the
roots $\epsi_1\pm\epsi_j$, $j=2,3,4$, so
$\sigma\alpha=-\epsi_1/2\pm\ldots<0$ for all $\alpha\in\Bu$. This
implies $\wt\Bu=\varnothing$. On the other hand, if $\alpha\in\Bu$,
then $\gamma=\beta-\alpha$ is not singular to any other root from
$D$ except $\beta$. Whence $f([x,e_{\gamma}])=\xi_{\beta}\cdot
x_{\alpha}\cdot N_{\alpha\gamma}\neq0$ if
$x=x_{\alpha}e_{\alpha}+\ldots\in\ut_{\Bu}$. Thus, $\bt=0$ as
required.

Thirdly, suppose there exists $j$ such that
$\beta=\epsi_1-\epsi_j\in D$ and $\epsi_1+\epsi_j\notin D$. In this
case, $\sigma\alpha>0$ if and only if
$\alpha=(\epsi_1+\epsi_j\pm\ldots)/2$. If $\gamma\in\wt\Phi^+$, then
$\alpha+\gamma\in\Bu$, so $f([e_{\alpha},e_{\gamma}])=0$. At the
same time if $\gamma\in\Bu$, then the coefficient of~$\epsi_j$ in
$\gamma$ is not less than $-1/2$, so $\alpha+\gamma\neq\beta$. On
the other hand, $\alpha+\gamma=\epsi_1\pm\ldots$, so
$\alpha+\gamma\notin S(\beta)$ for all $\beta\in D$. This yields
that $e_{\alpha}\in\bt$. But if $x\in\ut_{\Bu}$ and
$\alpha=(\epsi_1-\epsi_j\pm\ldots)/2\in\Supp{x}$, then
$\gamma=\beta-\alpha$ isn't singular to any other root from~$D$
except $\beta$. If $x=x_{\alpha}e_{\alpha}+\ldots$, then
$f([x,e_{\gamma}])=\xi_{\beta}\cdot x_{\alpha}\cdot
N_{\alpha\gamma}\neq0$, a~contradiction. Thus, $\bt=\langle
e_{\alpha},\alpha\in\wt\Bu\rangle_k$.

Similarly, if there exists $j$ such that $\beta=\epsi_1+\epsi_j\in
D$ and $\epsi_1-\epsi_j\notin D$, then $\sigma\alpha>0$ if and only
if $\alpha=(\epsi_1-\epsi_j\pm\ldots)/2$, i.e., $e_{\alpha}\in\bt$.
At the same time if $\alpha=(\epsi_1+\epsi_j\pm\ldots)$, then
$\gamma=\beta-\alpha$ isn't singular to any other root from~$D$
except~$\beta$, so $f([x,e_{\gamma}])\neq0$. It follows that if
$x\in\ut_{\Bu}$ and $\alpha\in\Supp{x}$, then $x\notin\at$. Hence
$\bt=\langle e_{\alpha},\alpha\in\wt\Bu\rangle_k$.

Finally, suppose $\epsi_1-\epsi_j,\epsi_1+\epsi_j\in D$ for some
$j$. Then $\sigma\alpha=-\epsi_1/2\pm\ldots<0$ for all
$\alpha\in\Bu$, so $\wt\Bu=\varnothing$. Let $\alpha$ be a root from
$\Bu$. Then $\alpha=(\epsi_1+z\cdot\epsi_j\pm\ldots)/2$, $z=\pm1$,
so $\alpha$ isn't singular to any other root from $D$ except
$\beta=\epsi_1+z\cdot\epsi_j$; this is also true for the root
$\gamma=\beta-\alpha$. Arguing as above, we see that $x\notin\at$ if
$x\in\ut_{\Bu}$ and $\alpha\in\Supp{x}$. Thus, $\bt=0$ as required.
The proof is complete.}

\medskip It follows from two previous Lemmas that
\begin{equation*}
\begin{split}
\dim\Omega&=\codim_{\ut}\at=\dim\ut-\dim\at=|\Phi^+|-(\dim\wt\at+\dim\bt)\\
&=|\wt\Phi^+|+|\Bu|-\dim\wt\at-\dim\bt=(|\wt\Phi^+|-\dim\wt\at)+|\Bu|-\dim\bt\\
&=(\dim\wt\ut-\dim\wt\at)+|\Bu|-\#\{\alpha\in\Bu\mid\sigma\alpha>0\}\\
&=\codim_{\wt\ut}\wt\at+\#\{\alpha\in\Bu\mid\sigma\alpha<0\}\\
&=\dim\wt\Omega+\Fo-\wt\Fo=\wt\Fo-\teta+\Fo-\wt\Fo=\Fo-\teta.
\end{split}
\end{equation*}
Therefore $\dim\Omega$ doesn't depend on~$\xi$ and is less or equal
to~$\Fo$. In other words, Theorem~\ref{mtheo} holds for all
orthogonal subsets of $\wt\Phi^+$.

\sst In this Subsection, we consider orthogonal subsets of $F_4^+$
which don't contain in $\wt\Phi^+$. In other words, we assume the
intersection of $D$ with $\Bu$ to be non-empty. It's easy to see
that if $\beta_1,\beta_2\in\Bu$ are orthogonal and $\beta_1\notin
S(\beta_2)$, then $\beta_2\in S(\beta_1)$, so without loss of
generality it can be assumed that $|D\cap\Bu|=1$ (see
Lemma~\ref{lemm_sing}). Clearly, $D$ doesn't contain the
roots~$\epsi_i$, $1\lee i\lee 4$; in other words, there exists a
unique short root contained in $D$.

The root system $F_4$ is self-dual, so there exists the bijection
$\vfi\colon F_4\to F_4$  such that $\vfi(F_4^+)=F_4^+$ and
$\vfi(S(\alpha))=S(\vfi(\alpha))$ for a given positive root
$\alpha$. Further,
$\vfi(\Bu)=\{\epsi_1\pm\epsi_j,{\epsi_2\pm\epsi_j}$, ${j=3,4}\}$
(signs are independent) and $\vfi(\wt\Phi)\cong C_4$. So if
$D\subset\vfi(\wt\Phi^+)$, then the results of~\cite{Ignatev3} can
be applied, and Theorem~\ref{mtheo} holds for the subset $D$.

Let
$\beta=(\epsi_1+z_2\cdot\epsi_2+z_3\cdot\epsi_3+z_4\cdot\epsi_4)/2\in
D$, $z_j=\pm1$. By previous remarks, $\wt D=D\setminus\{\beta\}$
coincides either with one of the subsets
$\{\epsi_1-z_3\cdot\epsi_3$, $\epsi_2-z_2z_4\cdot\epsi_4\}$,
$\{\epsi_1-z_4\cdot\epsi_3,\epsi_2-z_2z_3\cdot\epsi_3\}$ or with a
one-element subset of them. Since the root
$\beta=(\epsi_1-\epsi_2-\epsi_3-\epsi_4)/2$ is fundamental,
$\Omega=\xi_{\beta}e_{\beta}^*+\wt\Omega$, and so
$\dim\Omega=\dim\wt\Omega$, where $\wt\Omega=\Omega_{\wt
D,\wt\xi}\subset\ut^*$, $\wt\xi=\xi\setminus\{\xi_{\beta}\}$. On the
other hand, $r_{\beta}$ acts on the set of positive roots non-equal
to~$\beta$ by permutations, so
$l(\sigma)-s(\sigma)=(l(\wt\sigma)+1)-(s(\wt\sigma)+1)=l(\wt\sigma)-s(\wt\sigma)$,
where $\wt\sigma$ is the involution in $W$ corresponding to~$\wt D$.
But Theorem~\ref{mtheo} holds for the orbit~$\wt\Omega$, so we may
assume that $\beta\neq(\epsi_1-\epsi_2-\epsi_3-\epsi_4)/2$. By the
same argument, it can be assumed that the fundamental root
$\alpha_1=\epsi_2-\epsi_3$ doesn't belong to $D$ (if $|D|=2$, then
the problem reduces to elementary orbits; if $|D|=3$, then the
problem reduces to orbits associated with two-element subsets).

For a given $D$, denote by $\Mo\subset\Phi^+$ a subset satisfying
the following conditions. Firstly, if $\alpha+\gamma=\beta\in D$,
then $|\Mo\cap\{\alpha,\gamma\}|=1$. Secondly, for a given
$\gamma\in\Mo$, there exists $\alpha\in\Po$ such that
$\alpha+\gamma=\beta\in D$ (here $\Po=\Phi^+\setminus\Mo$). Thirdly,
$(\alpha+\Mo)\cap D$ consists either of the root~$\beta$ or of the
roots $\beta,\wt\beta=\alpha+\wt\gamma$, $\wt\gamma\in\Mo$, and in
the latter case $\wt\gamma\in S(\beta)$, $(\wt\alpha+\Mo)\cap
D=\{\beta\}$, where $\wt\alpha=\beta-\wt\gamma\in\Po$. Assume $\Mo$
exists. Then $\pt=\sum_{\alpha\in\Po}ke_{\alpha}$ is a maximal
isotropic subspace of the canonical form on $\Omega$, so
$\dim\Omega=2\cdot\codim_{\ut}\pt=2\cdot|\Mo|$ doesn't depend
on~$\xi$.

In the table below we list subsets $\Mo$ for all remaining $D$
(signs $\pm$ in the table are independent). It's straightforward to
check that they satisfy the above conditions. We also compute the
numbers $\Fo=l(\sigma)-s(\sigma)$ for all $D$. One can see that
$2\cdot|\Mo|\lee\Fo$ for all $D$. This concludes the proof of
Theorem~\ref{mtheo}.

\begin{center}
\begin{tabular}{||l|l|c|c||} \hline Subset $D$ &
Subset $\Mo$ & $|\Mo|$ &$\Fo$\\
\hline\hline 1) $\epsi_1+\epsi_3$, & $\epsi_1$, $\epsi_4$,
$\epsi_1-\epsi_2$, $\epsi_1\pm\epsi_4$, & 7 & 14\\
$(\epsi_1-\epsi_2-\epsi_3+\epsi_4)/2$ &
$(\epsi_1+\epsi_2+\epsi_3\pm\epsi_4)/2$ & &\\
\hline 2) $\epsi_1-\epsi_4$, & $\epsi_1-\epsi_2$, $\epsi_1-\epsi_3$,
& 4& 8\\
$(\epsi_1-\epsi_2-\epsi_3+\epsi_4)/2$ &
$(\epsi_1\pm\epsi_2-\epsi_3-\epsi_4)/2$ & &\\
\hline 3) $\epsi_2+\epsi_4$, & $\epsi_4$, $\epsi_2-\epsi_3$ & 2 &
4\\$(\epsi_1-\epsi_2-\epsi_3+\epsi_4)/2$ & & &\\
\hline 4) $\epsi_1-\epsi_3$, & $\epsi_1-\epsi_2$, & 3&  6\\
$(\epsi_1-\epsi_2+\epsi_3-\epsi_4)/2$ &
$(\epsi_1-\epsi_2-\epsi_3\pm\epsi_4)/2$ & &\\
\hline 5) $\epsi_1+\epsi_4$, & $\epsi_1$, $\epsi_3$,
$\epsi_1-\epsi_2$, $\epsi_1-\epsi_3$, & 6& 12\\
$(\epsi_1-\epsi_2+\epsi_3-\epsi_4)/2$ &
$(\epsi_1-\epsi_2\pm\epsi_3+\epsi_4)/2$ & &\\
\hline 6) $\epsi_2+\epsi_3$, &
$\epsi_3$, $\epsi_3\pm\epsi_4$ & 3& 6\\
$(\epsi_1-\epsi_2+\epsi_3-\epsi_4)/2$ & & &\\
\hline 7) $\epsi_2-\epsi_4$, & $\epsi_3-\epsi_4$, & 2& 4\\
$(\epsi_1-\epsi_2+\epsi_3-\epsi_4)/2$ &
$(\epsi_1-\epsi_2-\epsi_3-\epsi_4)/2$ & &\\
\hline 8) $\epsi_1-\epsi_3$, & $\epsi_4$, $\epsi_1-\epsi_2$, & 4& 8\\
$(\epsi_1-\epsi_2+\epsi_3+\epsi_4)/2$ &
$(\epsi_1-\epsi_2-\epsi_3\pm\epsi_4)/2$ & &\\
\hline 9) $\epsi_1-\epsi_4$, & $\epsi_3$, $\epsi_1-\epsi_2$,
$\epsi_1-\epsi_3$, & 5& 10\\
$(\epsi_1-\epsi_2+\epsi_3+\epsi_4)/2$ &
$(\epsi_1-\epsi_2\pm\epsi_3-\epsi_4)/2$ & &\\
\hline 10) $\epsi_2+\epsi_3$, & $\epsi_3$, $\epsi_3\pm\epsi_4$, & 4& 8\\
$(\epsi_1-\epsi_2+\epsi_3+\epsi_4)/2$ &
$(\epsi_1-\epsi_2+\epsi_3-\epsi_4)/2$ & &\\
\hline 11) $\epsi_2+\epsi_4$, & $\epsi_4$, $\epsi_3+\epsi_4$, & 3&
6\\
$(\epsi_1-\epsi_2+\epsi_3+\epsi_4)/2$ &
$(\epsi_1-\epsi_2-\epsi_3+\epsi_4)/2$ & &\\
\hline 12) $\epsi_1+\epsi_3$, &
$\epsi_1$, $\epsi_2+\epsi_3$, $\epsi_1\pm\epsi_4$, & 8& 16\\
$(\epsi_1+\epsi_2-\epsi_3-\epsi_4)/2$ &
$(\epsi_1-\epsi_2\pm\epsi_3\pm\epsi_4)/2$ & &\\
\hline
\end{tabular}
\end{center}

\begin{center}
\begin{tabular}{||l|l|c|c||} \hline Subset $D$ &
Subset $\Mo$ & $|\Mo|$ &$\Fo$\\
\hline\hline 13) $\epsi_1+\epsi_4$, & $\epsi_1$, $\epsi_1-\epsi_2$,
$\epsi_1-\epsi_3$, & 7& 14\\
$(\epsi_1+\epsi_2-\epsi_3-\epsi_4)/2$ &
$(\epsi_1-\epsi_2\pm\epsi_3\pm\epsi_4)/2$ & &\\
\hline 14) $\epsi_2+\epsi_3$, &
$\epsi_2$, $\epsi_2-\epsi_3$, $\epsi_2\pm\epsi_4$ & 4& 8\\
$(\epsi_1+\epsi_2-\epsi_3-\epsi_4)/2$ & & &\\
\hline 15) $\epsi_2+\epsi_4$, & $\epsi_2$, $\epsi_2-\epsi_3$,
$\epsi_2-\epsi_4$ &3& 6\\
$(\epsi_1+\epsi_2-\epsi_3-\epsi_4)/2$ & & &\\
\hline 16) $\epsi_1+\epsi_3$, &
$\epsi_1$, $\epsi_4$, $\epsi_1-\epsi_2$, $\epsi_1\pm\epsi_4$, &9& 18\\
$(\epsi_1+\epsi_2-\epsi_3+\epsi_4)/2$ &
$(\epsi_1-\epsi_2\pm\epsi_3\pm\epsi_4)/2$ & &\\
\hline 17) $\epsi_1-\epsi_4$, & $\epsi_2$, $\epsi_1-\epsi_2$, $\epsi_1-\epsi_3$,& 6& 12\\
$(\epsi_1+\epsi_2-\epsi_3+\epsi_4)/2$ & $\epsi_2-\epsi_3$, $(\epsi_1\pm\epsi_2-\epsi_3-\epsi_4)/2$& &\\
\hline 18) $\epsi_2+\epsi_3$, & $\epsi_2$, $\epsi_2\pm\epsi_4$, & 5& 10\\
$(\epsi_1+\epsi_2-\epsi_3+\epsi_4)/2$ & $\epsi_4$,
$\epsi_2-\epsi_3$ & &\\
\hline 19) $\epsi_2-\epsi_4$, & $\epsi_2$, $\epsi_4$,
$\epsi_2-\epsi_3$, $\epsi_2+\epsi_4$ &4& 8\\
$(\epsi_1+\epsi_2-\epsi_3+\epsi_4)/2$ & & &\\
\hline 20) $\epsi_1-\epsi_3$, & $\epsi_2$, $\epsi_2-\epsi_3$, $\epsi_2-\epsi_4$, &6 & 14\\
$(\epsi_1+\epsi_2+\epsi_3-\epsi_4)/2$ & $\epsi_3-\epsi_4$, $(\epsi_1\pm\epsi_2-\epsi_3-\epsi_4)/2$ & &\\
\hline 21) $\epsi_1+\epsi_4$, & $\epsi_1$, $\epsi_2$,
$\epsi_1-\epsi_3$, $\epsi_2\pm\epsi_4$, &9 & 20\\
$(\epsi_1+\epsi_2+\epsi_3-\epsi_4)/2$ &
$(\epsi_1\pm\epsi_2-\epsi_3\pm\epsi_4)/2$ & &\\
\hline 22) $\epsi_2+\epsi_4$, & $\epsi_2$, $\epsi_3$, $\epsi_2\pm\epsi_3$, & 6& 12\\
$(\epsi_1+\epsi_2+\epsi_3-\epsi_4)/2$ & $\epsi_2-\epsi_4$,
$\epsi_3-\epsi_4$ & &\\
\hline 23) $\epsi_1-\epsi_3$, & $\epsi_2$, $\epsi_4$,
$\epsi_2\pm\epsi_3$, $\epsi_2+\epsi_4$, &7 &16\\
$(\epsi_1+\epsi_2+\epsi_3+\epsi_4)/2$ &
$(\epsi_1+\epsi_2-\epsi_3\pm\epsi_4)/2$ & &\\
\hline 24) $\epsi_1-\epsi_4$, &
$\epsi_2$, $\epsi_3$, $\epsi_2+\epsi_3$, $\epsi_3-\epsi_4$, & 8& 18\\
$(\epsi_1+\epsi_2+\epsi_3+\epsi_4)/2$ & $\epsi_2\pm\epsi_4$, $(\epsi_1+\epsi_2\pm\epsi_3-\epsi_4)/2$ & &\\
\hline 25) $\epsi_2-\epsi_4$, & $\epsi_2$, $\epsi_2+\epsi_3$, $\epsi_2+\epsi_4$, & 7& 14\\
$(\epsi_1+\epsi_2+\epsi_3+\epsi_4)/2$ & $\epsi_3+\epsi_4$, $\epsi_3$, $\epsi_4$, $\epsi_2-\epsi_3$ & &\\
\hline 26) $\epsi_1-\epsi_3$, $\epsi_2-\epsi_4$, & $\epsi_2$,
$\epsi_4$,
$\epsi_2+\epsi_4$, $\epsi_2\pm\epsi_3$, & 7& 16\\
$(\epsi_1+\epsi_2+\epsi_3+\epsi_4)/2$ &
$(\epsi_1+\epsi_2-\epsi_3\pm\epsi_4)/2$ & &\\
\hline 27) $\epsi_1-\epsi_3$, $\epsi_2+\epsi_4$, & $\epsi_2$,
$\epsi_2-\epsi_3$, $\epsi_2-\epsi_4$, & 6& 14\\
$(\epsi_1+\epsi_2+\epsi_3-\epsi_4)/2$ & $\epsi_3-\epsi_4$, $(\epsi_1\pm\epsi_2-\epsi_3-\epsi_4)/2$ & &\\
\hline 28) $\epsi_1+\epsi_3$, $\epsi_2+\epsi_4$, & $\epsi_1$,
$\epsi_4$, $\epsi_1-\epsi_2$, $\epsi_3\pm\epsi_4$, & 7& 14\\
$(\epsi_1-\epsi_2-\epsi_3+\epsi_4)/2$ &
$(\epsi_1-\epsi_2+\epsi_3\pm\epsi_4)/2$ & &\\
\hline 29) $\epsi_1-\epsi_3$, $\epsi_2+\epsi_4$, & $\epsi_2$,
$\epsi_2-\epsi_3$, & 4& 10\\
$(\epsi_1-\epsi_2+\epsi_3-\epsi_4)/2$ &
$(\epsi_1-\epsi_2-\epsi_3\pm\epsi_4)/2$ & &\\
\hline 30) $\epsi_1+\epsi_4$, $\epsi_2+\epsi_3$, & $\epsi_3$,
$\epsi_4$,
$\epsi_2+\epsi_4$, $\epsi_3+\epsi_4$, & 6& 14\\
$(\epsi_1-\epsi_2+\epsi_3-\epsi_4)/2$ &
$(\epsi_1-\epsi_2\pm\epsi_3+\epsi_4)/2$ & &\\
\hline 31) $\epsi_1-\epsi_3$, $\epsi_2+\epsi_4$, & $\epsi_4$,
$\epsi_2-\epsi_3$, & 4& 10\\
$(\epsi_1-\epsi_2+\epsi_3+\epsi_4)/2$ &
$(\epsi_1-\epsi_2-\epsi_3\pm\epsi_4)/2$ & &\\
\hline 32) $\epsi_1-\epsi_4$, $\epsi_2+\epsi_3$, & $\epsi_3$,
$\epsi_2-\epsi_4$, $\epsi_3-\epsi_4$, & 5& 12\\
$(\epsi_1-\epsi_2+\epsi_3+\epsi_4)/2$ &
$(\epsi_1-\epsi_2\pm\epsi_3-\epsi_4)/2$ & &\\
\hline 33) $\epsi_1+\epsi_3$, $\epsi_2+\epsi_4$, & $\epsi_1$,
$\epsi_2$, $\epsi_2-\epsi_4$, $\epsi_3\pm\epsi_4$, & 8& 18\\
$(\epsi_1+\epsi_2-\epsi_3-\epsi_4)/2$ & $\epsi_2+\epsi_3$, $(\epsi_1-\epsi_2+\epsi_3\pm\epsi_4)/2$ & &\\
\hline 34) $\epsi_1+\epsi_4$, $\epsi_2+\epsi_3$, &
$\epsi_1$, $\epsi_2$, $\epsi_2-\epsi_3$, $\epsi_3+\epsi_4$, & 7& 16\\
$(\epsi_1+\epsi_2-\epsi_3-\epsi_4)/2$ & $\epsi_2+\epsi_4$, $(\epsi_1-\epsi_2\pm\epsi_3+\epsi_4)/2$ & &\\
\hline 35) $\epsi_1+\epsi_3$, $\epsi_2-\epsi_4$, & $\epsi_1$,
$\epsi_4$, $\epsi_1-\epsi_2$, $\epsi_3\pm\epsi_4$, & 9& 20\\
$(\epsi_1+\epsi_2-\epsi_3+\epsi_4)/2$ &
$(\epsi_1-\epsi_2\pm\epsi_3\pm\epsi_4)/2$ & &\\
\hline 36) $\epsi_1-\epsi_4$, $\epsi_2+\epsi_3$, &
$\epsi_2$, $\epsi_2-\epsi_3$, $\epsi_2-\epsi_4$, & 6& 14\\
$(\epsi_1+\epsi_2-\epsi_3+\epsi_4)/2$ &
$\epsi_3-\epsi_4$, $(\epsi_1\pm\epsi_2-\epsi_3-\epsi_4)/2$ & &\\
\hline
\end{tabular}
\end{center}

\end{document}